\documentclass{article} 
\usepackage[a4paper, total={6in, 8in}]{geometry}  
\usepackage{latexsym}
\usepackage{graphicx}
\usepackage{mathptmx}
\usepackage[T1]{fontenc}

\usepackage{amsmath}
\usepackage{amsfonts}
\usepackage{amssymb}
\usepackage{amsbsy}
\usepackage{amsthm}

\usepackage{caption}
\usepackage{subcaption} 
\usepackage{comment}

\newcommand*{\email}[1]{\href{mailto:#1}{\nolinkurl{#1}} } 
\usepackage[pdftex,colorlinks=true,urlcolor=blue,citecolor=black,anchorcolor=black,linkcolor=black]{hyperref}

\newtheoremstyle{wsc}
{3pt}
{3pt}
{}
{}
{\bf}
{}
{.5em}
{}

\theoremstyle{wsc}
\newtheorem{theorem}{Theorem}

\newtheorem{lemma}[theorem]{Lemma}
\usepackage[toc,page]{appendix}
\begin{document}
\date{} 
%
%

\title{Patient Assignment and Prioritization for Multi-Stage Care with Reentrance}

\author{
Wei Liu$^1$~~~
Mengshi Lu$^1$~~~
Pengyi Shi$^1$~~~
\smallskip 
\\
$^1$Mitchell E. Daniels, Jr. School of Business, Purdue University, West Lafayette, Indiana, USA
}

\maketitle

\vspace{-12pt}

\section*{ABSTRACT}  
In this paper, we study a queueing model that incorporates patient reentrance to reflect patients' recurring requests for nurse care and their rest periods between these requests. Within this framework, we address two levels of decision-making: the priority discipline decision for each nurse and the nurse-patient assignment problem. We introduce the shortest-first and longest-first rules in the priority discipline decision problem and show the condition under which each policy excels through theoretical analysis and comprehensive simulations. For the nurse-patient assignment problem, we propose two heuristic policies. We show that the policy maximizing the immediate decrease in holding costs outperforms the alternative policy, which considers the long-term aggregate holding cost. Additionally, both proposed policies significantly surpass the benchmark policy, which does not utilize queue length information. 

\section{INTRODUCTION} 
\label{sec:intro}

Nurse staffing and workload management are critical to efficient hospital operations and timely, high-quality patient care delivery. The COVID-19 pandemic not only exacerbated the shortage of nurses and resources, but also shed light on the issue of physical beds versus staffed beds. In the US, many nurses left their hospital employment for traveling contracts due to the increased and unbalanced workload. Nursing shortage became significantly amplified due to nurse turnover during the COVID-19 pandemic. The implications of inadequate staffing are far-reaching, including nurse burnout, high turnover rates, and a decline in the quality of care. A key factor that drives nurse burnout is the unbalanced and volatile workload they face during their shifts. Nurse burnout leads to higher turnover, which then exacerbates the nursing shortage and understaffing, which results in a vicious cycle that negatively affects the well-being of healthcare professionals and compromises patient outcomes. Motivated by the burgeoning needs of nurse workload management, we use queueing-based simulation methods in this paper to study balancing nurse workload efficiently while enhancing patient care and operational cost-effectiveness.

Traditional nurse staffing models focus on the patient count (number of patients) to measure the nurse workload, which does not account for the fact that patients often need help several times during their stay and patients are highly heterogeneous in terms of the number of requests: more severe patients may have higher medical needs and generate more frequent requests with a longer length-of-stay (LOS) in the hospital. A larger proportion of high-need patients could significantly increase a nurse’s workload, even if the number of patients attended by this nurse may be similar to another nurse who manages mostly low-need patients. In other words, traditional models that focus on patient count could underestimate how busy nurses are because they do not account for the repeated care requests from patients and the heterogeneity in these repeated patient requisitions.

In this paper, we introduce an improved model to measure nursing workloads beyond simply counting the number of patients. We build a queueing model with returns that captures patients' repeated requests for nurses and their time spent in bed between these requests. In other words, a patient may be fine one moment and require immediate attention the next, then settle down again once their needs are met. By factoring in these recurring care requests, our model presents a more accurate picture of what nurses face during their shifts. Based on this model, we study two specific decision problems. First, how should nurses prioritize which patients to serve next (e.g., those who are just admitted or those who are close to being discharged)? Second, under the different prioritization rules that each nurse may use, when a new patient shows up to the unit, who should this patient be assigned to so that a more balanced workload among nurses can be achieved (typically a decision to be made by the charge nurse of a unit). By studying these two questions, we hope to also generate insights into how hospitals can better determine the right number of nurses needed, helping to prevent nurse burnout and ensuring patients receive the timely care they require. 

Our research is most related to the research on queueing models with reentry. Specifically, both \cite{yom2014erlang} and \cite{van2016restricted} study the scenarios where patients can reenter the needy queue after receiving care from a nurse. However, they assume that all patients are of a single type and that each patient will leave the system with a fixed probability after receiving service from a nurse. In contrast, our study examines a more realistic setting by considering multiple patient types, each potentially requiring varying amounts of service from the nurse. Moreover, our research focuses on the coordination of nurse-patient assignments and the prioritization of patients for each nurse, markedly diverging from the main focus of \cite{yom2014erlang} and \cite{van2016restricted} on nurse staffing levels. This differentiation distinctly sets our research apart and contributes unique insights to the field.

This paper is organized as follows. 
Section 2 introduces the model setting and the two-level decision structure (the priority discipline decision for each nurse and the nurse-patient assignment problem). Section~3 studies the priority discipline decision problem. Section 4 studies the nurse-patient assignment problem. Section 5 shows simulation results where we compare various policies for both levels of decisions through an extensive simulation. Section 6 concludes the paper. The complete description of notations used in this paper is in Appendix \ref{app:holdingcost} and all proofs are provided in Appendix \ref{app:proof}.

\section{The Model} 
\label{sec:model}

\begin{figure}[htb] 
\centering
\includegraphics[width=0.6\textwidth]{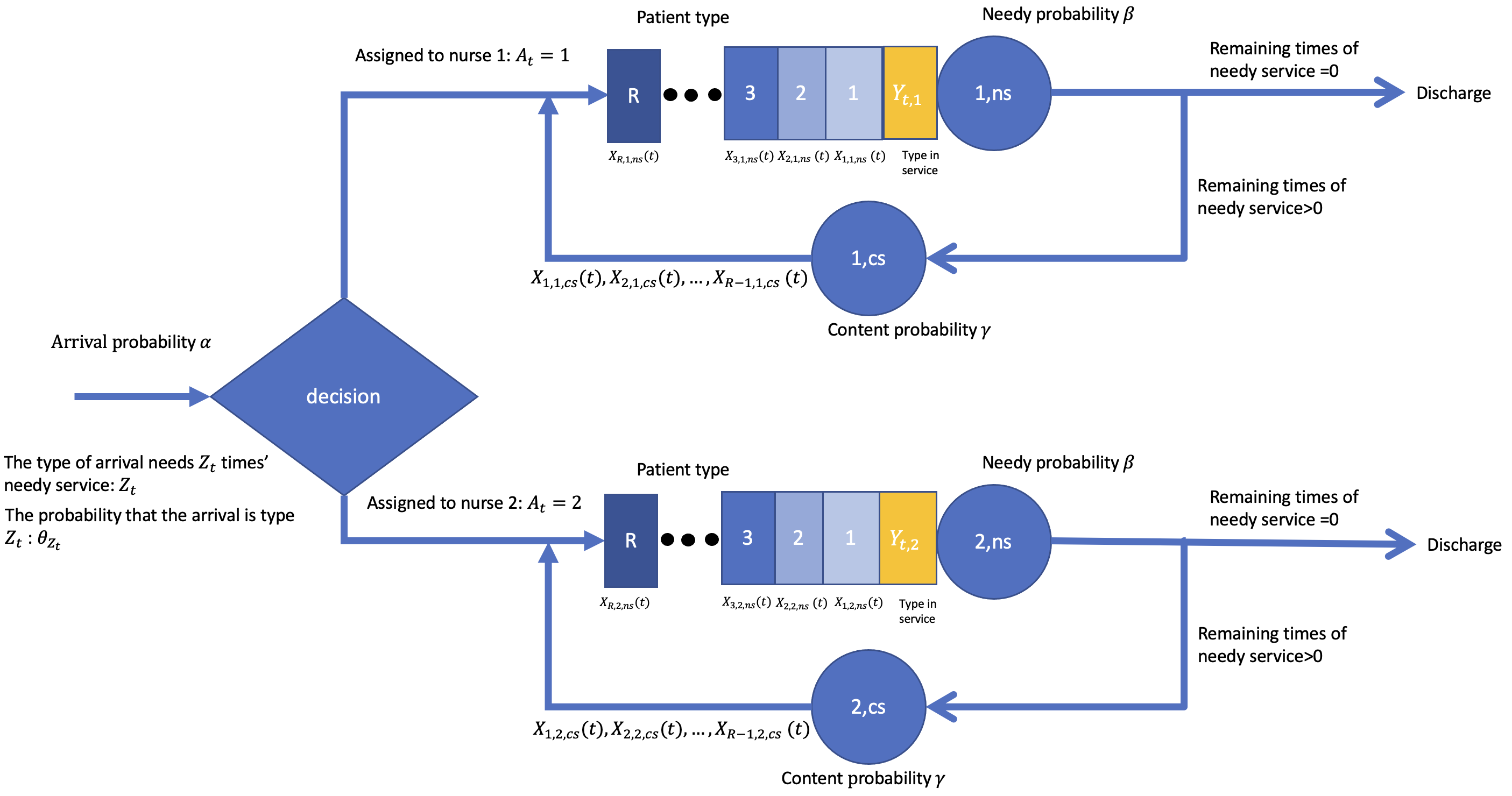}
\caption{Diagram for the two levels of decision processes: nurse-patient assignment decision (between the two nurses) and priority decision (among patients of one nurse).}
\label{fig:flow-chart}
\end{figure}

We consider a nurse-patient assignment queuing system with two levels of decisions. The first level is for an incoming patient, which nurse to assign this patient to. This decision is usually made by the unit charge nurse. The second level is for each nurse in the unit, among all existing patients, which patient she should prioritize to proceed. 
The diagram is illustrated in Figure \ref{fig:flow-chart} with two unit nurses as an example.

\noindent\textbf{System dynamics and states. }
We consider a discrete-time system setup with predetermined periods $1\le t\le T$. Patients arrive in a discrete-time manner. At the beginning of each period $t$, one or zero patient arrives in this period. 
The probability of a patient's arrival is denoted by $\alpha$ ($0\le \alpha\le 1$). The arrival patient could be one of $R$ types. We use $r \in \{1, 2, \dots, R\}$ to denote the patient type, which means that the patient needs $r$ times (``stages'') of service from a nurse. We use $\theta_{r}$ to denote the probability that a new patient is of type $r$, and use the random variable $Z_t$ to denote the type of arrival in period $t$ if one arrival occurs. 

An incoming patient is first assigned to a unit nurse and joins the queue of the assigned nurse, awaiting the nurse's care. We refer to such a status as the \emph{needy} state when patients are waiting to receive nurse care, following the convention in the literature \cite{yom2014erlang,van2016restricted}.
Let $\beta$ be the probability that the nurse can finish processing the service needs of a patient in one discrete period. 
To ensure the stability of the system, we need to make sure the inflow rate is less than the outflow rate, that is, $\alpha < \frac{\beta I }{\sum_{r=1}^R r\theta_r}$, where $I$ is the number of nurses, and $\sum_{r=1}^R r\theta_r$ is the expected number of times to visit the nurse. 

Once the nurse takes care of the patient's current-stage needs, the patient leaves the nurse queue. 
The patient changes into the \emph{content} state and their stage changes from $r$ to $r-1$ if $r>1$; otherwise, the patient is discharged from the hospital when all service stages are finished. 
A content patient stays in bed, {with $\gamma>0$ as the probability to generate the next stage's service request.}

To track patients in different stages of their service in each nurse's queue, we separate the patients in the needy state versus those in the content state. For patients who are assigned to nurse $i$ and are in the needy state (``ns''), we use $X_{r,i,ns} (t) \in \mathbb{N}_0 $ to denote the number of patients with $r\in[R]$ stages' left right after the patient arrival at the beginning of period $t$, where  $\mathbb{N}_0 $ is the set of nonnegative integers.  
Let ${\bf X}_{i,ns}  (t) = [X_{1,i,ns} (t), X_{2,i,ns}  (t),\cdots, X_{R,i,ns}  (t) ]$ be the vector for these patient counts in the needy state for nurse $i\in [I]$. 
Similarly, for patients who are assigned to nurse $i$ and are currently in the content state (``cs''), we use $X_{r,i,cs}  (t)  \in \mathbb{N}_0$ to denote the number of patients with $r\in[R-1]$ stages left. Let ${\bf X}_{i,cs}  (t)  = [X_{1,i,cs}  (t) , X_{2,i,cs}  (t),\cdots, X_{R-1,i,cs}  (t) ]$ be the vector for these patient counts in the content state for nurse $i$. We note that patients in the content state have a maximum of $R-1$ remaining stages instead of $R$. This comes from our assumption that new arrivals always go through the nurse queue first (corresponding to admission to hospitals). Thus, when the patients enter the content state, they have already had one service interaction with the nurse, which reduces their maximum possible remaining stages to $R-1$. 
The notations are summarized in Table \ref{tab:notation} in Appendix \ref{app:holdingcost}.

\noindent\textbf{Two levels of decisions. }
Upon a patient arrival, the charge nurse determines which nurse to assign this patient to. This is the first level of decision. 
We let $A_t$ be the decision of the charge nurse. In the two-nurse example as illustrated in Figure~\ref{fig:flow-chart}, $A_t=1$ means that we assign the patient to nurse 1, while $A_t=2$ means that the patient is assigned to nurse 2. We reserve $A_t=0$ for the situation when there is no arrival in period $t$ and no assignment decision needs to be made.

Beyond the initial nurse-patient assignment determined by the charge nurse, each unit nurse needs to prioritize among various patient types to determine who to serve next. This is the second level of decision. 
That is, each nurse needs to determine $Y_{t,i} \in [R]$,  the type of patients for service in period $t$ in needy state at nurse $i$ (with $Y_{t,i}  =0$ when no patient in service). We use ${\bf Y}_{t}   =  [Y_{t,1}  , Y_{t,2}  , \cdots, Y_{t,I}  ] $ for the vector of service decisions across all nurses. We assume a non-preemptive priority queue in the needy state.

\begin{figure}[htb]
\centering
\includegraphics[width=0.6\textwidth]{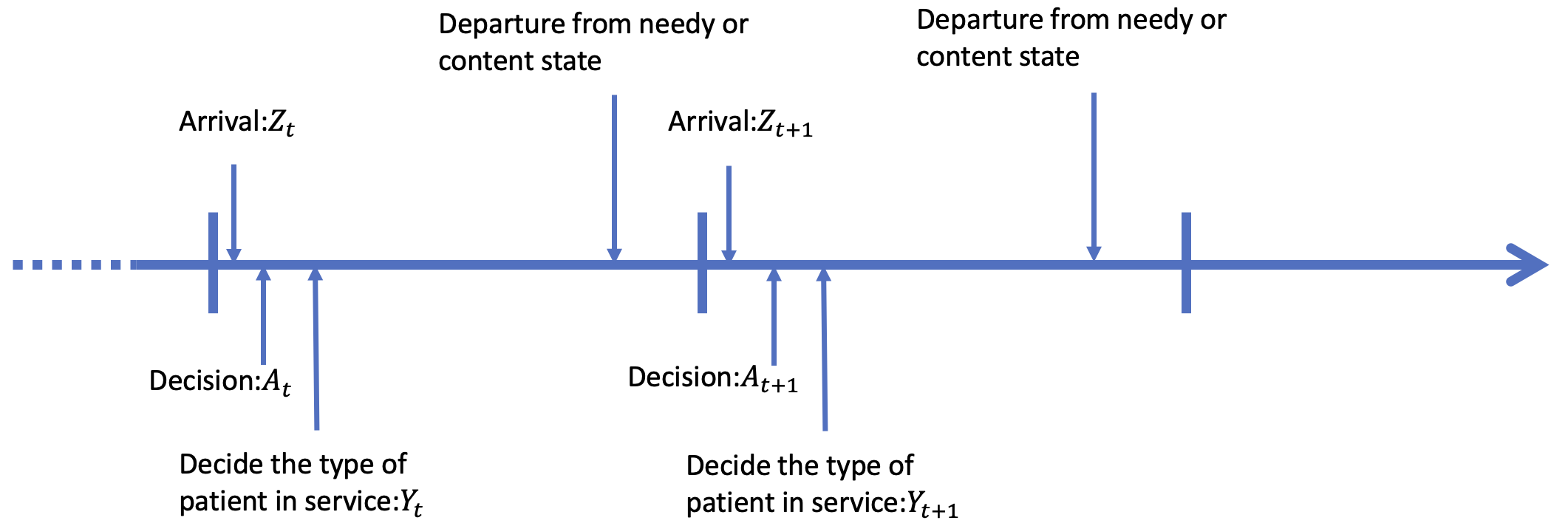}
\caption{Timeline for the nurse-patient assignment problem and priority discipline decision problem}
\label{fig:timeline}
\end{figure}

\noindent\textbf{Decision timeline. } 
We illustrate the timeline of the decision-making process in Figure~\ref{fig:timeline}. 
At the start of period $t$, a patient of type $Z_t$ arrives and joins the queue. Following this, the charge nurse executes decision $A_t$, assigning this incoming patient to a specific nurse. 
Subsequently, each of the nurses determines the patient type $Y_{t,i} $ they will serve in this epoch.  
Note that we assume the new patient who just arrived and was assigned to a nurse in period $t$ can receive the service for the assigned nurse in period $t$. At the end of period $t$, the patient currently being served by the nurse may finish service with probability $\beta$, who would then either enter the content state and stay in the bed, or get discharged if all service stages are finished. With probability $1-\beta$, the patient continues getting service from the nurse. Each patient currently in the content state may enter the needy state (requesting service from the nurse) with probability $\gamma$.  

To solve this hierarchical two-level decision problem, we proceed in two steps. In the first step, we decide the type of patient for service, namely, the priority discipline decision problem, and then in the second step, we solve the nurse-patient assignment problem. We show the details in Section \ref{sec:priority} and  \ref{sec:np_assign}, respectively. 

\section{Priority discipline decision problem}  \label{sec:priority}
In this section, we focus on each individual nurse and study the priority decision, i.e., how each nurse should decide which type of patients to serve. Let $\Pi_i$ be the set of feasible policies for nurse $i$. We assume that the holding cost in one period at the needy state for one type $r$ patient is $h(r,a)= r^a$ for some parameter $a\in[0,\infty]$; there is no holding cost in the content state. We can interpret the holding cost as the potential harm to patients, such as increased treatment costs due to delays incurred when they wait for nurse service. Different choices of $a$ enable us to model a variety of scenarios. For example, when $a=0$, the holding cost is constant for patients in any stage; as $a$ increases, the model increasingly differentiates patients by the remaining stages in their hospital stay, with higher-type (potentially more critical) patients incurring exponentially higher holding costs. {The parameter $a$ plays a crucial role in allowing us to capture potential nonlinear effects brought by delay during different stages of patient care (e.g., the delay during the admission stage has a more significant impact than the delay during the discharge stage) and the heterogeneous effect of delay among different patients, which are salient features relevant in healthcare settings \cite{chan2012optimizing}. 
In practice, estimating the holding cost parameter $a$ is difficult. Hospital managers can choose the value of $a$ by considering the tradeoff among various performance metrics, such as the queue length for all patients versus the queue length for patients with more extensive remaining service needs. We provide details of choosing $a$ according to such a tradeoff in Appendix~\ref{app:holdingcost}.
}


The total expected cost under policy $\pi \in \Pi_i $ for nurse $i$ is  
$$ J_i^\pi = E \bigg[ \sum_{t=1}^T\sum_{r=1}^R h(r,a)   { X}^\pi_{r,i,ns}  (t) \bigg] ,$$  
where $ { X}^\pi_{r, i,ns} (t) $ is the number of type $r$ patients under policy $\pi$ in needy state at nurse $i$ at the beginning of period $t$.  
Let $ J_i^{\pi^*} = \min_{\pi \in \Pi } J_i^\pi $ be the optimal total expected  cost. 
Note that the priority discipline decision problem is a finite-horizon discrete-time Markov decision process. 
Solving this problem is challenging due to the large state space and we focus on static priority policies. To address this challenge, we consider two policies, $\pi^1$ and $\pi^2$, where $\pi^1$ is the \textit{shortest-first policy}, prioritizing patients based on the shortest remaining service times; and $\pi^2$ is the \textit{longest-first policy}, prioritizing patients with the longest remaining service times.  
These two policies are proposed based on the following Lemma.  
  \begin{lemma} \label{lemma:convexity}  
For any given policy $\pi$ that is independent of $a$, ${J_i^{\pi}}$ is increasing and convex in $a$. 
\end{lemma}   

From Lemma \ref{lemma:convexity}, we see that with the increase of $a$,  the cost associated with patients having longer remaining service times in the needy queue increases. 
This highlights the substantial impact of the value of $a$ on the optimal policy.  
Consequently, this suggests a strategic approach: employing the shortest-first rule might be more effective when $a$ is small, whereas the longest-first rule could be preferable as $a$ increases.  
We show the optimality of these two policies in the next subsection.

\subsection{Two-patient queueing system}  
To derive analytical properties, we consider a two-patient queueing system. In this simplified system, there are only two patients, one of type $i$ and the other of type $j$  ($i \leq j$), both just arrived (starting from the needy state). Furthermore, we assume that (i) there are no external arrivals and there are enough periods left to finish the service of these two patients, i.e., consider a \textit{clearing system}; (ii) the service times for either the needy or content states are deterministic in each stage and are the same for both patients.

Without loss of generality, let the type $i$ patient be referred to as patient 1, and the type $j$ patient as patient 2. We define system $s_1$ as the system that serves patient 1 first, followed by serving patient 2. Similarly, we define system $s_2$ as the one that serves patient 2 first and then serves patient 1. 
It is important to note that no additional decisions need to be made following this initial selection.
Essentially, system $s_1$ employs the shortest-first rule, while system $s_2$ employs the longest-first rule.

Let $S_{m, 1,  ns, l }$ be the deterministic service time of patient $1$ in needy state when she joins the nurse queue for the $l$th time in system $s_m$, {with $l=i$ being the last service this patient needs from the nurse}; similarly $S_{m, 2,  ns, l }$ for patient 2, {with $l=j$ being the last service this patient needs from the nurse}. 
Further, let $D_{m, i,  ns, l }$ be the departure time when patient $i$ leaves the needy state for the $l$th time in system $s_m$ for $i\in \{1, 2\}$.    


\begin{lemma} \label{lemma:discharge_sequence}
For a two-patient queueing system, where patient 1 is of type $i$ and patient 2 is of type $j$  ($i \leq j$), we have 
\begin{equation}
\begin{aligned}
 D_{1,  1,  ns, i  }      & \le D_{2,  1,  ns, i   }      \le   D_{1,  2,  ns, j   }   , \\ 
 D_{1,  1,  ns, i  }    & \le D_{2,  2,  ns, j   }       \le  D_{1,  2,  ns, j   }    . 
\end{aligned} 
\end{equation}
\end{lemma}   

We see that the shortest-first rule (system $s_1$) tends to make the patient with fewer times of remaining service leave first, and increases the length of stay for the patient with longer remaining times of service (even longer than the other patient when she is switched to system $s_2$). It implies that under the shortest-first rule, the patient with longer remaining times of service tends to linger in the nurse queue (the needy state). On the other hand, the longest-first rule tends to increase the length of stay for the patient with shorter remaining times of service and decrease the length of stay for the patient with longer remaining times of service. It implies that under the longest-first rule, the patient with fewer remaining times of service tends to linger in the nurse queue (the needy state). However, since the per-period holding cost differs by the remaining service stages, it is unclear which rule incurs a higher total holding cost. In the following theorem, we compare the total cost in systems $s_1$ and $s_2$, generating insights into the comparison over the shortest-first and longest-first rules.   
\begin{theorem}  \label{theorem:threshold} 
There exists a threshold of $a=\hat a \in [0,1]$ such that the total cost incurred in $s_1$ is smaller than that under $s_2$ when $a\le \hat a$, and the total cost incurred in $s_2$ is smaller than that under $s_1$ when $a> \hat a$. 
\end{theorem}   

Although our theoretical result is for the simplified two-patient system, it generates insights for more general settings, which we will demonstrate using simulation. That is, under the shortest-first rule, the patients with longer remaining times of service tend to accumulate in the nurse queue, whereas under the longest-first rule, the patients with less remaining times of service tend to accumulate in the nurse queue. Due to this difference, the number of patients under the longest-first rule would also be higher than that under the shortest-first rule (since patients are not prioritized to be discharged). As a result, when $a$ is small, e.g., when $a=0$ with the per-period holding cost being the same for patients in any stage, it is better to use the shortest-first rule since the total queue length under the longest-first rule is larger. However, as $a$ increases, it becomes more and more costly for the patients with longer remaining times of service to stay in the queue, e.g., when $a=1$, the per-period holding cost increases in the remaining stages $r$. Once $a$ goes beyond a threshold, the additional holding cost incurred by the additional patients with longer remaining times of service under the shortest-first rule outweighs the additional holding cost incurred by the longer queue length under the longest-first rule. Consequently, the longest-first rule becomes better.   

\section{Nurse-patient assignment problem}  \label{sec:np_assign}

Having studied the priority problem, we now move to the next level of decision and study the nurse-patient assignment problem. Let $\Pi^{np}$ be the set of feasible policies for the charge nurse. The total expected cost under policy $\pi^{np} \in \Pi^{np} $ is    
$$ J^{\pi^{np}}  = E \bigg[ \sum_{i=1}^I   \sum_{t=1}^T\sum_{r=1}^R h(r,a)   { X}^{\pi^{np}}_{r, i,ns} (t) \bigg] ,$$   
where $ { X}^{\pi^{np}}_{r, i,ns} (t) $  is the number of type $r$ patients under policy $\pi^{np} $ in needy state at nurse $i$ at the beginning of period $t$.  
Let $ J^{\pi^{np*}} = \min_{\pi^{np} \in \Pi } J^{\pi^{np}} $ be the optimal total expected  cost. 
This is a finite-horizon discrete-time Markov decision process as well.

Given the large state space, traditional solution methods become impractical. We introduce and compare two heuristic approaches for this decision problem. These two heuristics are motivated by the $c\mu $ rule and maximum pressure policy~\cite{mandelbaum2004scheduling,dai2005maximum}.  
Specifically, we define Heuristic 1 (H1) as follows: at the start of each time period, we calculate the following for each nurse \(i \in [I]\): 
$$ c_i^{np,1} =   \sum_{r=1}^R h(r,a)   { X}_{r, t,i,ns}. $$  
Here, the patient is allocated to the nurse with the lowest \(c_i^{np,1}\) value. In cases where multiple nurses have identical \(c_i^{np,1}\) values, the patient is assigned to any one of these nurses at random, with each nurse having an equal probability of being selected. The primary objective of Heuristic 1 is to identify the nurse who can, from a myopic perspective, lead to a larger decrease in the instantaneous holding cost.  
The state $X_{r, t,i,ns}$ provides a snapshot of the system's current status. This focus on maximizing the instantaneous cost reduction is beneficial, particularly during periods of high congestion. 


In Heuristic 2 (H2), at the beginning of each period, we compute 
$$ c_i^{np,2} =   \sum_{r=1}^R   ( { X}_{r, t,i,ns}  + { X}_{r, t,i,cs}  )    \sum_{r'=1}^r   h(r', a)   $$ 
for nurse $i\in [I]$. 
The patient is assigned to the nurse with the lowest \(c_i^{np,2}\). In instances where \(c_i^{np,2}\) values are equal for some nurses, the patient allocation is randomized among these nurses with equal probability for each. 
In contrast, Heuristic 2 looks into the future by taking into account the aggregate holding cost associated with the needy service for all patients, whether they are in needy or content states.  
This heuristic is based on the fact that patient of type \(r\) will be serviced by nurses for \(r\) times, thereby adopting a more ``long-term view'' to evaluating holding costs. Comparing Heuristic 2 with Heuristic 1, it goes beyond the myopic view but may underestimate the immediate need to reduce the holding cost as Heuristic 1 does. 

To demonstrate the efficacy of the two aforementioned heuristics, we establish a baseline policy wherein each incoming patient is assigned randomly to one of the two nurses with equal probability. This approach does not leverage any available information regarding the patient counts associated with each nurse. By implementing this baseline, we aim to showcase the better performance of our proposed heuristics. 
We show the comparison mainly via simulation study; see section \ref{subsec:sim_np}.


%

\section{Simulation}  
In this section, we compare the performance of the shortest-first and the longest-first patient prioritization rules, as well as two heuristics for the nurse-patient assignment problem.
\subsection{Patient Prioritization}
\label{subsec:sim_pp} 

To examine the efficacy of the patient prioritization heuristic strategies, in this subsection, we assume that all patient arrivals are directly assigned to a single nurse. Patients are classified into five categories; that is, $R=5$. The probabilities of each patient category are set as $[\theta_{1}, \theta_{2}, \theta_{3},\theta_{4}, \theta_{5}]= [0, 0.3380,0.2238, 0.1481, 0.0981]$. These probabilities are derived by normalizing the geometric distribution used in the numerical experiments as reported in \cite{yom2014erlang}.
There are no category 1 arrivals (i.e., the arrival probability equals zero) because all patients need at least two nurse services.
For each simulation, we run 20 replications, with 10,000 periods in each replication and the first 2000 periods being excluded. The cost under each policy is computed based on the average of these 20 replications \cite{law2007simulation}.

\begin{figure}[htb]
\centering
\includegraphics[width=0.3\textwidth]{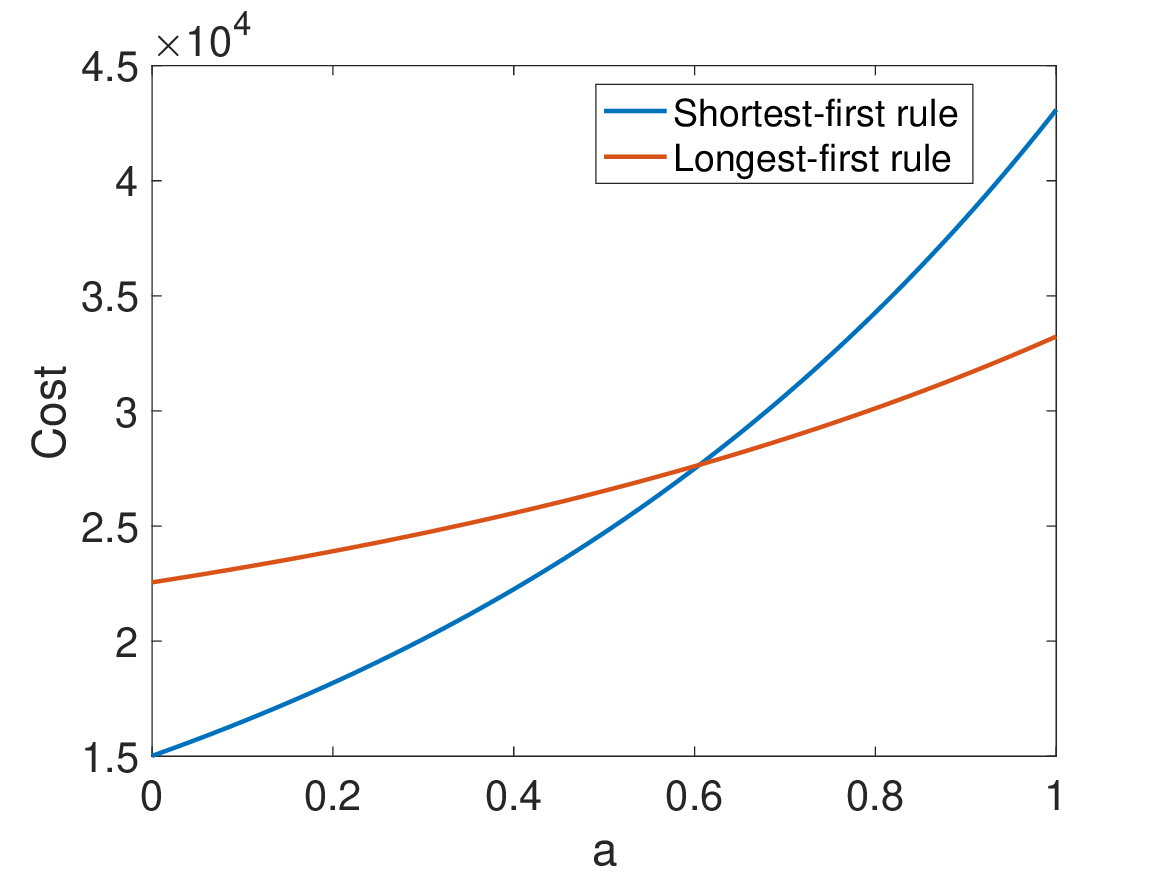}
\caption{Average total holding costs using the two priority rules impacted by $a$} 
\label{fig_cost_a_short_long} 
\end{figure}

First, we investigate the impact of the holding cost parameter $a$.
Recall that the holding cost of a type $r$ patient is equal to $r^a$.
When $a=0$, all patients have the same holding cost.
When $a>0$, patients with more remaining nurse visits incur higher holding costs.
The cost difference is more significant with a larger $a$.
Figure \ref{fig_cost_a_short_long} shows the total holding costs for the two priority rules as a function of $a$, with $\alpha=0.2$, $\beta=0.8$, and $\gamma=0.1$. 
We see that the holding costs associated with both heuristics exhibit convexity with respect to $a$. 
Notably, there exists a threshold $\hat a = 0.61$, when $a < \hat a$, the shortest-first rule results in lower costs, and when $a > \hat a$, the longest-first rule performs better.
These findings align with Lemma \ref{lemma:convexity} and Theorem \ref{theorem:threshold}.

\begin{figure}
     \centering
     \begin{subfigure}[b]{0.3\textwidth}
         \centering
         \includegraphics[width=\textwidth]{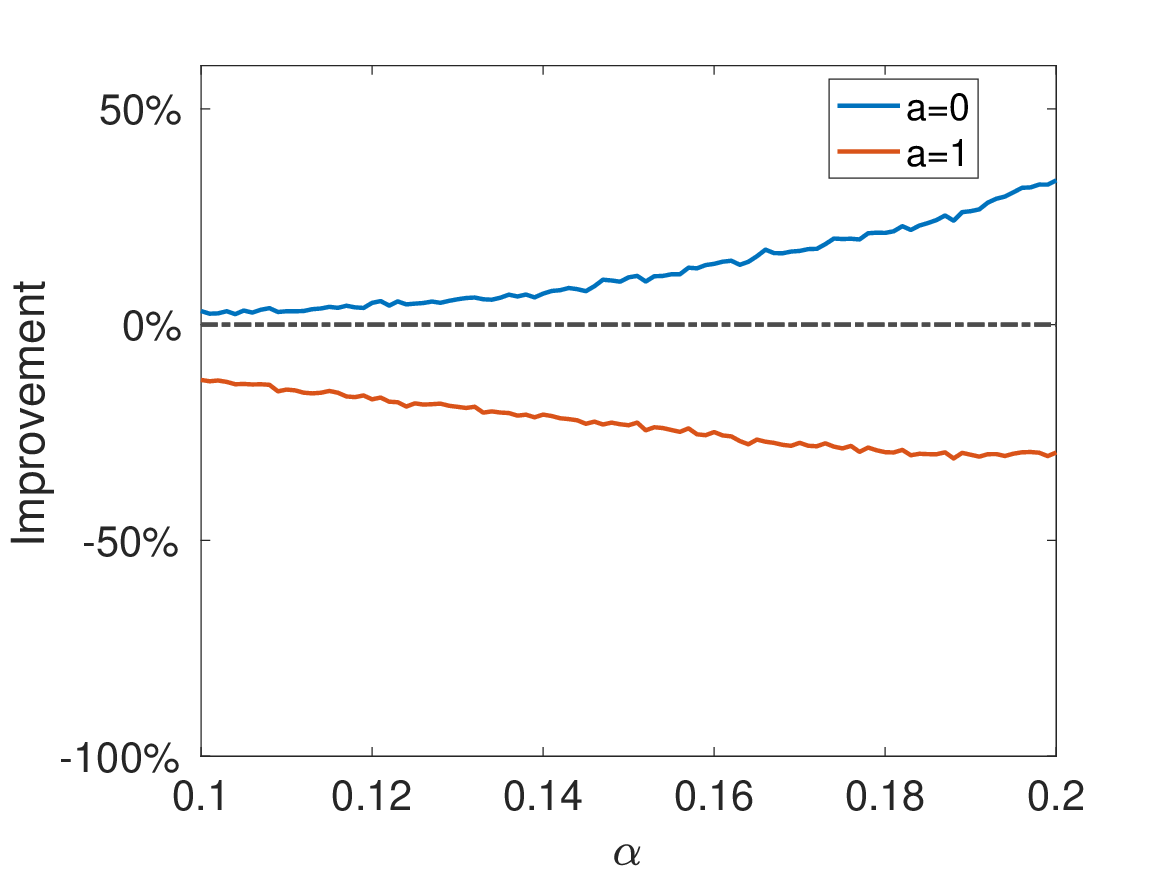}
         \caption{Improvement impacted by $\alpha$ }
         \label{fig_costinc_alpha_a01_beta02_gamma01_short_long}
     \end{subfigure}
     \hfill
     \begin{subfigure}[b]{0.3\textwidth}
         \centering
         \includegraphics[width=\textwidth]{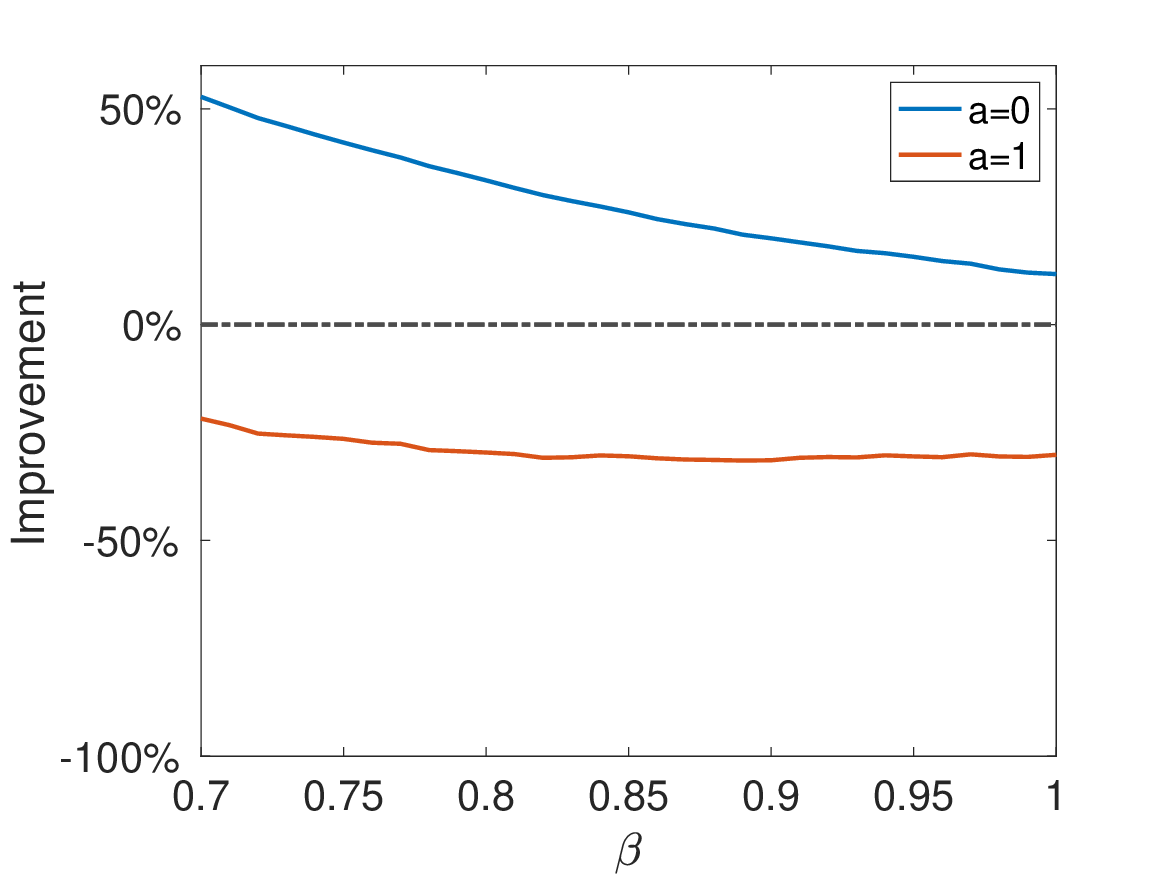}
         \caption{Improvement impacted by $\beta$}   
         \label{fig_costinc_beta_a01_alpha02_gamma01_short_long}
     \end{subfigure}
          \hfill
     \begin{subfigure}[b]{0.3\textwidth}
         \centering
         \includegraphics[width=\textwidth]{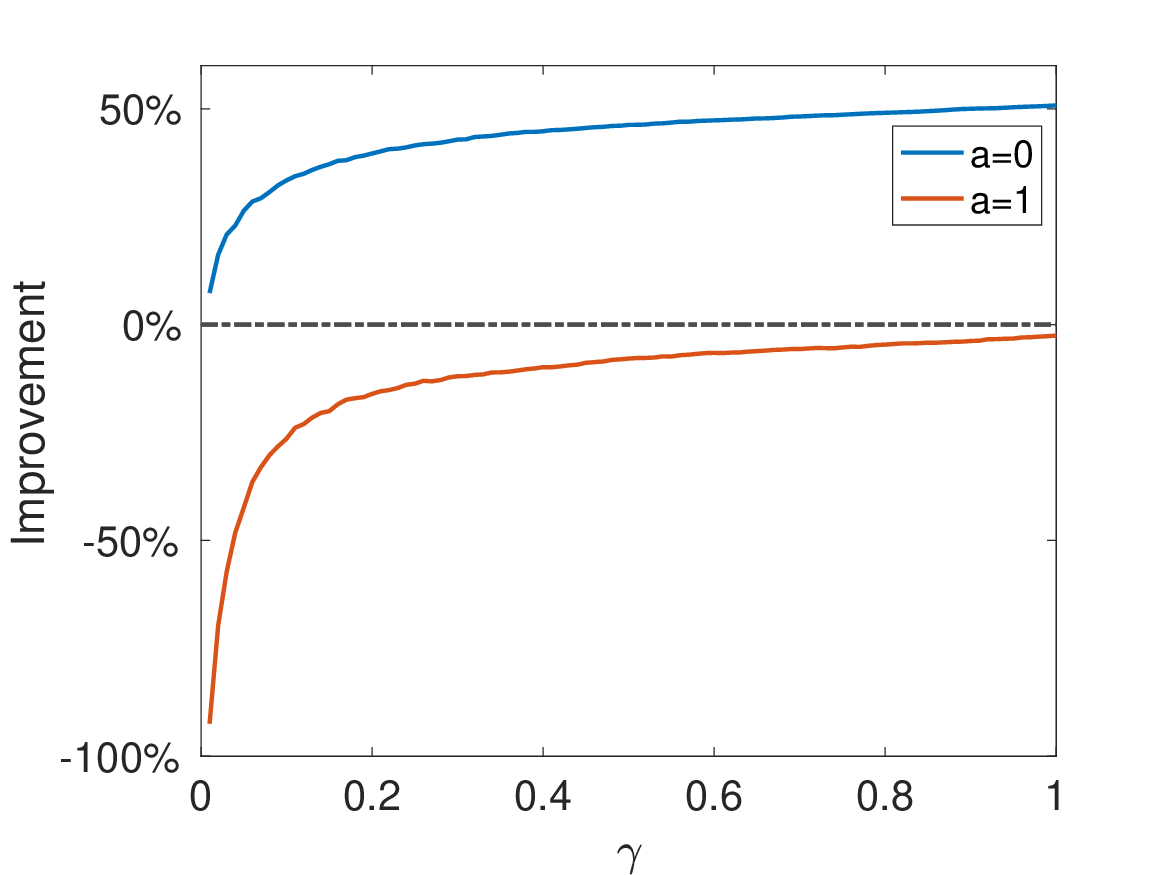}
         \caption{Improvement impacted by $\gamma$}   
         \label{fig_costinc_gamma_a01_beta02_alpha02_short_long}
     \end{subfigure}
        \caption{Improvement of the shortest-first rule over the longest-first rule}
        \label{fig_np_alpha_a01_gamma01_beta08}  
\end{figure}

Next, we study how the comparison between the shortest-first and longest-first rules is impacted by parameters $\alpha$, $\beta$, and $\gamma$.
We focus on two cases of holding costs with $a=0$ or 1 and measure the relative improvement (i.e., percentage of cost reduction) of the shortest-first rule over the longest-first rule.
As shown in Theorem \ref{theorem:threshold} and Figure \ref{fig_cost_a_short_long}, compared to the longest-first rule, the shortest-first rule performs better when $a=0$ and worse when $a=1$.
Therefore, when $a=0$, the improvement (cost reduction) is positive; when $a=1$, the ``improvement'' is negative, and its absolute value is the percentage of cost increase of the shortest-first rule over the longest-first rule. 
Figure \ref{fig_costinc_alpha_a01_beta02_gamma01_short_long} illustrates the improvement of the shortest-first rule as a function of $\alpha$, with $\beta = 0.8$ and $\gamma=0.1$.
We see that when $a=0$, the improvement tends to increase in $\alpha$.   
This is because, with more arrivals, the congestion under the longest-first rule is becoming increasingly worse than that under the shortest-first rule, resulting in higher holding costs.
When $a=1$, the cost increase (i.e., the absolute value of the negative improvement) of the shortest-first rule also tends to increase in $\alpha$.
This is because when $a=1$, patients with longer remaining service time incur higher holding costs. With higher arrival rates, the impact of the higher holding costs is further amplified, causing the shortest-first rule to perform increasingly worse.
Similarly, Figure \ref{fig_costinc_beta_a01_alpha02_gamma01_short_long} shows the improvement as a function of $\beta$, with $\alpha=0.2$ and $\gamma=0.1$, and Figure \ref{fig_costinc_gamma_a01_beta02_alpha02_short_long} illustrates the improvement as a function of $\gamma$, with $\alpha=0.2$ and $\beta=0.8$. 
The corresponding monotonicity can be similarly explained by the change in congestion levels that results from variations in $\beta$ and $\gamma$.


\subsection{Nurse-Patient Assignment}    \label{subsec:sim_np}  

In this subsection, we compare the two nurse-patient assignment heuristic policies (H1 and H2) with the baseline policy that randomly assigns patients to each nurse with equal probabilities.
We consider the case with two nurses. Patients assigned to the same nurse are served according to the shortest-first or longest-first priority rules. Similar to the previous subsection, we focus on two holding cost settings with ($a=0$ or 1). Recall that the shortest-first patient prioritization rule performs better when $a=0$ and the longest-first rule performs better when $a=1$. Therefore, we apply the shortest-first rule when $a=0$ and the longest-first rule when $a=1$. We adhere to the settings outlined in Section \ref{subsec:sim_pp} for patient categories and simulation replications.


We first show how the improvements of the heuristics H1 and H2 over the baseline policy depend on the arrival probability $\alpha$, with $\gamma=0.1$ and $\beta = 0.8$.
As shown in Figure \ref{fig_np_alpha_a01_gamma01_beta08}, regardless of the value of $a$, improvements of the heuristics generally increase in $\alpha$.
This is because as $\alpha$ increases, there is more congestion in the system, which highlights the effect of patient assignment.
Moreover, we see that H1 performs better than H2, and the difference is also increasing in $\alpha$. 
It implies that while H2 accounts for long-term holding costs, H1 provides a more immediate reflection of the system's state. This real-time perspective is particularly crucial in scenarios of high congestion, where the ability to respond swiftly to changing conditions can significantly enhance performance. 

\begin{figure}
     \centering
     \begin{subfigure}[b]{0.3\textwidth}
         \centering
         \includegraphics[width=\textwidth]{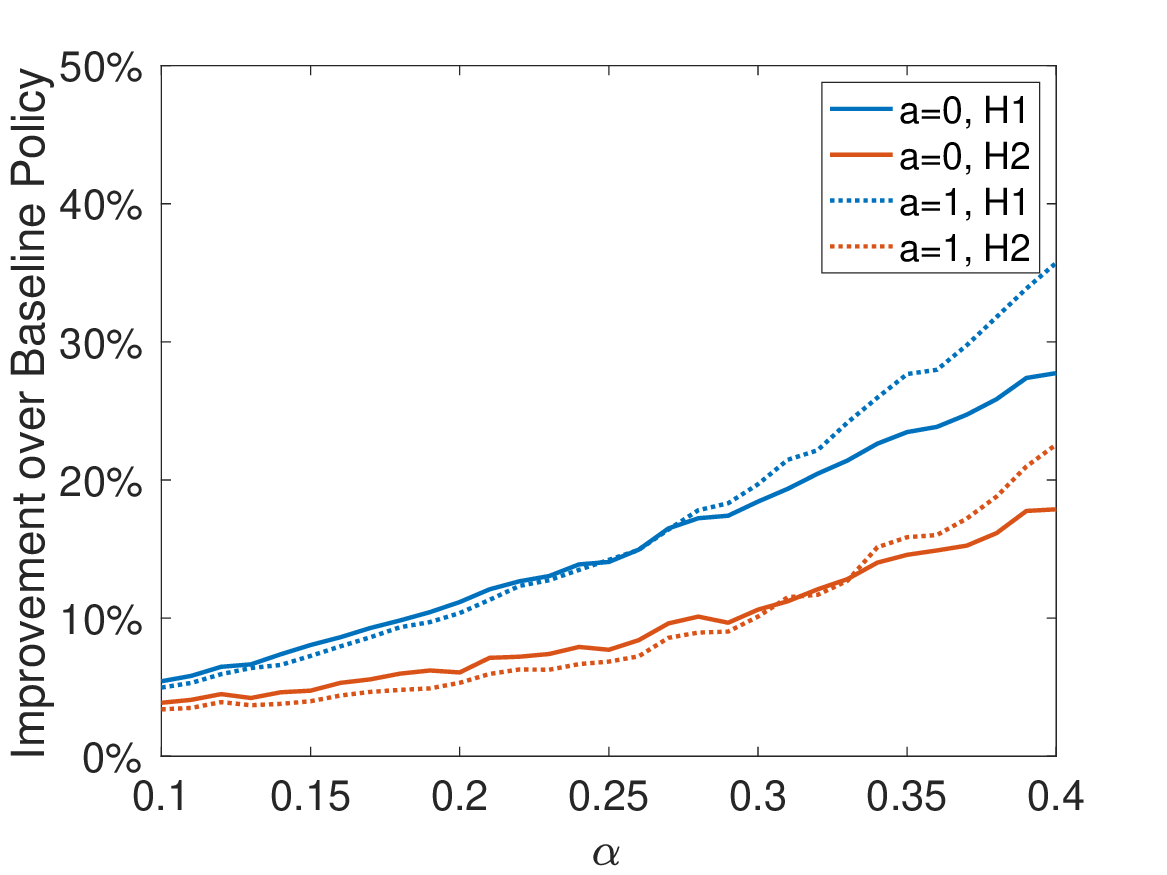}
         \caption{Improvement impacted by $\alpha$ }
         \label{fig_np_alpha_a01_gamma01_beta08}
     \end{subfigure}
     \hfill
     \begin{subfigure}[b]{0.3\textwidth}
         \centering
         \includegraphics[width=\textwidth]{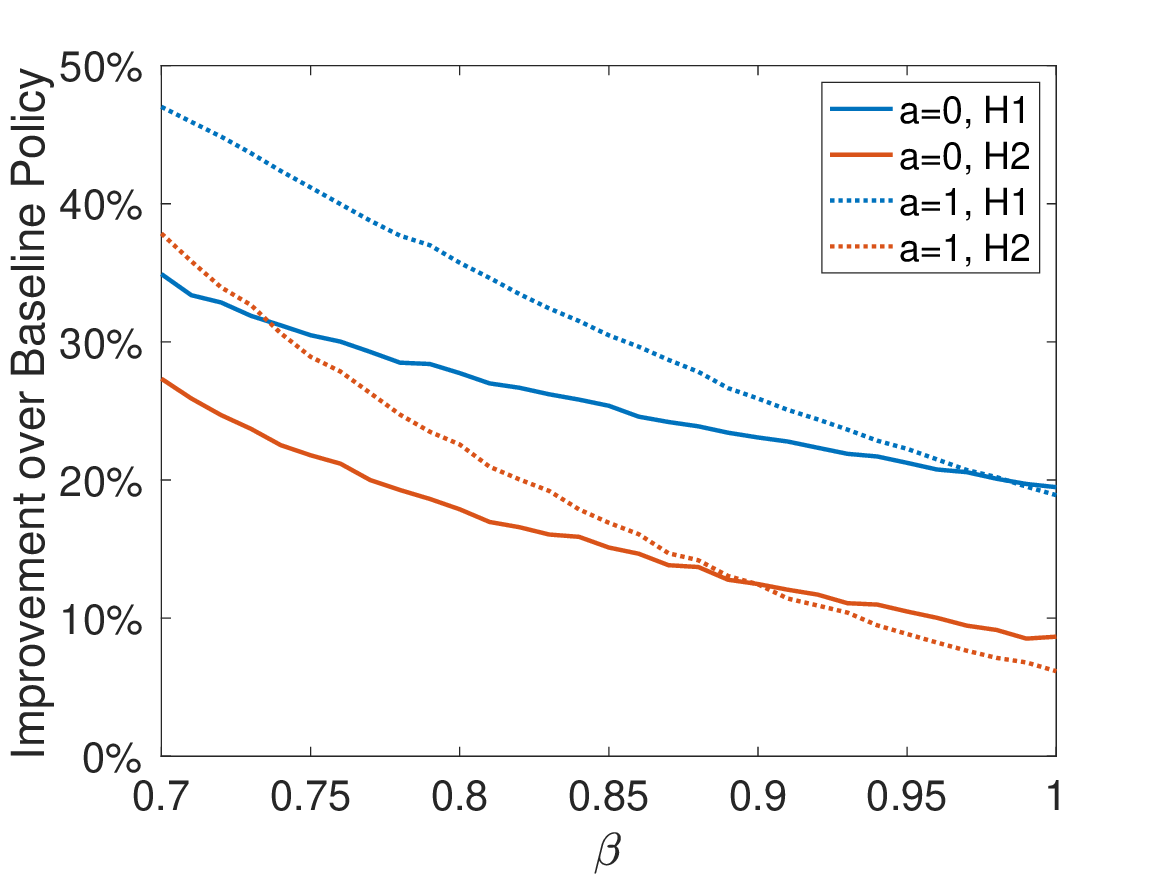}
         \caption{Improvement impacted by $\beta$}   
         \label{fig_np_beta_a01_gamma01_alpha04}
     \end{subfigure}
          \hfill
     \begin{subfigure}[b]{0.3\textwidth}
         \centering
         \includegraphics[width=\textwidth]{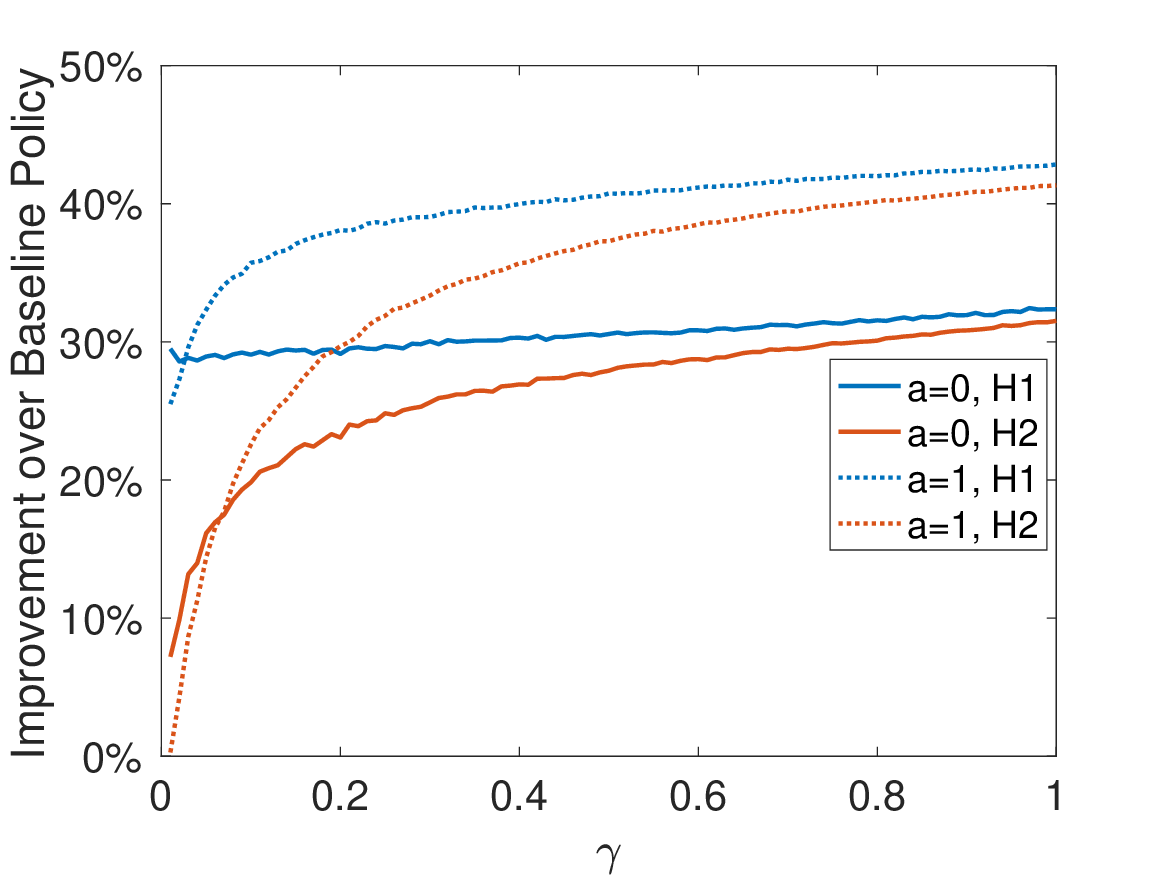}
         \caption{Improvement impacted by $\gamma$}   
         \label{fig_np_gamma_a01_alpha04_beta08}
     \end{subfigure}
        \caption{Improvement of the heuristics over random assignment}
        \label{fig_np_alpha_beta_gamma_a01}  
\end{figure}

Similarly, we study the impact of $\beta$, with $\gamma=0.1$ and $\alpha=0.4$, in Figure \ref{fig_np_beta_a01_gamma01_alpha04}, and show the impact of $\gamma$, with $\alpha=0.4$ and $\beta=0.8$ in Figure \ref{fig_np_gamma_a01_alpha04_beta08}.  
We see that as $\beta$ increases or $\gamma$ decreases, the improvements of the heuristics decrease.
This is consistent with our previous observation that the heuristics lead to more improvement when the system is more congested.
Also similar to our previous observation, H1 still outperforms H2.
Moreover, the difference between H1 and H2 widens as $\beta$ increases or as $\gamma$ decreases. 
This is because when $\beta$ increases or $\gamma$ decreases, there are more patients in the content state, which amplifies the difference between the two heuristics.

\section{Conclusions and Future Work}  

In this paper, we analyze a nurse-patient assignment queuing system characterized by two distinct decision-making levels. The first level concerns nurse-patient assignments, wherein the charge nurse determines the most suitable nurse for each incoming patient. The second level addresses the priority discipline decision problem, focusing on the priority in which a nurse attends to her assigned patients. 
To address this hierarchical problem, our approach unfolds in two phases. Initially, we tackle the priority discipline decision problem, selecting the patient type to prioritize. Subsequently, we address the nurse-patient assignment challenge.
Within the priority discipline framework, we evaluate the efficacy of both shortest-first and longest-first policies, delineating the scenarios under which each policy is better. For the nurse-patient assignment dimension, we introduce two strategic policies. Our analysis demonstrates that the policy emphasizing immediate holding costs outperforms its counterpart by offering a nuanced and timely snapshot of the system's state. Such an acute understanding is invaluable, especially under high congestion, where rapid response to fluctuating conditions is essential to optimize system performance. 
{ 
Future research directions include exploring scenarios where the charge nurse has only partial information (not knowing exactly) about the number of service requests when patients enter the queue. Additionally, collecting real operational data and calibrating our model with it, beyond using the current synthetic data simulation, would be valuable.}

\bibliographystyle{plain} 
\bibliography{literature}
\section*{AUTHOR BIOGRAPHIES}


\noindent {\bf \MakeUppercase{Wei Liu}} is a Ross-Lynn Postdoctoral Fellow in the Daniels School of Business at Purdue University. His research interests include data-driven resource allocation in real world operations, including airline, healthcare, and service operations, to enhance network resilience and service level.
His email address is \email{liu3568@purdue.edu} and his website is \url{https://sites.google.com/view/wei-liu-stor/home}. \\

\noindent {\bf \MakeUppercase{Mengshi Lu}} is an Associate Professor in the Daniels School of Business at Purdue University. His research interests include supply chain risk management, operations of innovative service systems, and project management. His email address is \email{mengshilu@purdue.edu} and his website is \url{https://business.purdue.edu/directory/bio.php?username=lu420}. \\

\noindent {\bf \MakeUppercase{Pengyi Shi}} is an Associate Professor of Operations Management in the Daniels School of Business at Purdue University. Her
research focuses on building data-driven, high-fidelity models and developing predictive and prescriptive analytics to support decision-making under uncertainty in healthcare and the criminal justice system. Her email is \email{shi178@purdue.edu} and her website is \url{https://web.ics.purdue.edu/~shi178}.

\appendix 

\begin{table}[ht]
\begin{minipage}[b]{0.5\linewidth}
  \centering
  \caption{Description of Model Notations}
\scalebox{0.75}{
    \begin{tabular}{lp{30em}} 
    \hline
    Notations & \multicolumn{1}{l}{Description} \\ 
    \hline
    $t$   & Period \\
    $T$   & Total periods in consideration \\
        $I$   & Number of nurses \\
    $r$   & Patient type \\
    $R$   & Number of patient types \\
    $\theta_r$ & Probability that a new arrival patient is of type $r$ \\
    $Z_t$ & Type of arrival in period $t$ \\
    $\alpha$ & Arrival probability in each period \\
    $\beta$ & Nurse service completion probability in each period \\
    $\gamma$ & Content-to-needy probability in each period \\
    $X_{r,i,ns}  (t) $ & Patient counts in needy state with $r$ stages' left at nurse $i$ in period $t$ \\
    $X_{r,i,cs}  (t) $ & Patient counts in content state with $r$ stages left at nurse $i$ in period $t$ \\
    $A_t$ & Decision of the charge nurse in period $t$ \\
    $Y_{t,i} $ & Patient type for service in needy state at nurse $i$  in period~$t$ \\
    \hline 
    \end{tabular}%
} 
  \label{tab:notation}%
\end{minipage}\hfill
\begin{minipage}[b]{0.35\linewidth}
\centering
\includegraphics[width=0.95\textwidth]{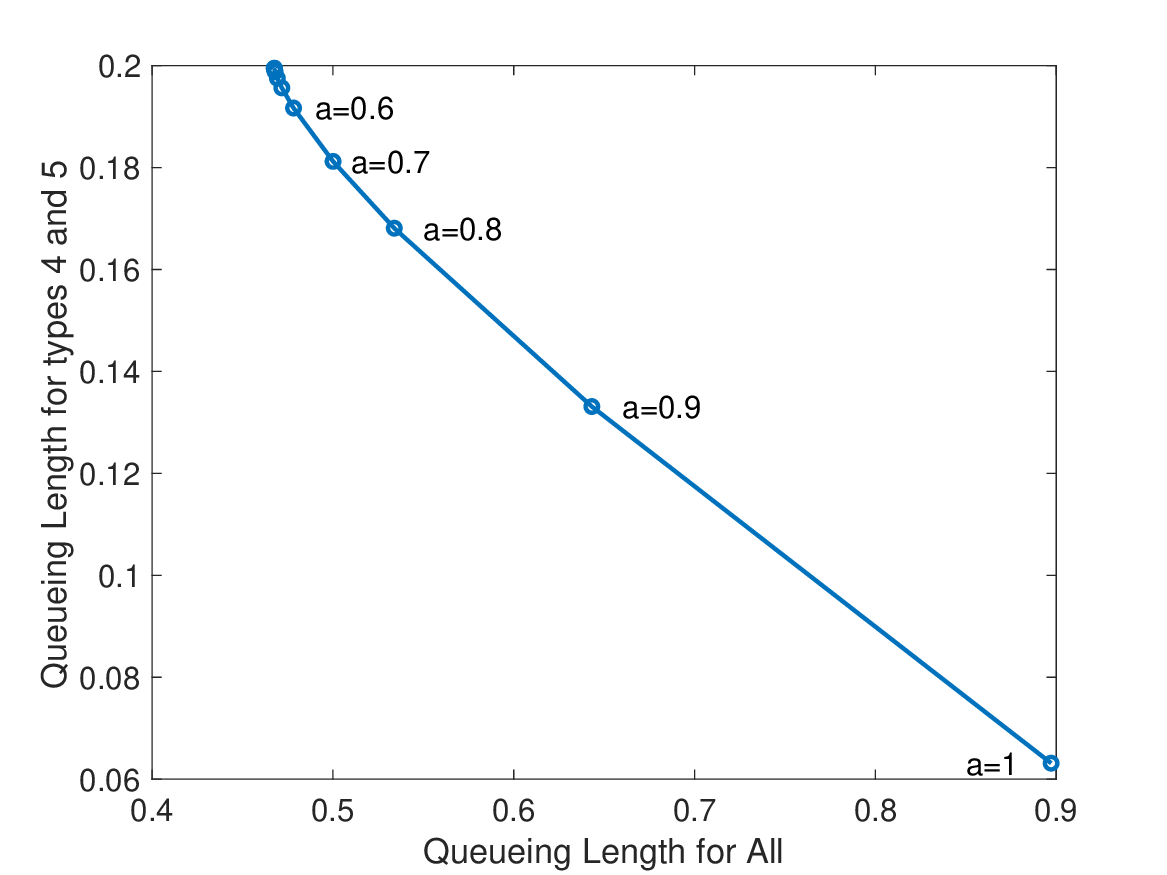}
\captionof{figure}{\small{Average Queue Length for All Patients Against Average Queue Length for Types 4 and 5 Patients} } 
\label{fig:parameter_tune} 
\end{minipage}
\end{table}

\section{Notations and Holding cost parameter selection} \label{app:holdingcost} 
{
Table~\ref{tab:notation} summarizes the model notations. Next, we demonstrate how to choose the parameter $a$ in the holding cost function $h(r, a)$ using tradeoff curves. We consider $R=5$ as in Section 5, and we analyze the tradeoff between the queue length for all patients and the queue length for more severe patients (types 4 and 5). When $a$ is zero, the unit holding cost is the same across all patients. As $a$ increases, priority is given to those with larger remaining stages (types 4 and 5) due to higher holding costs. Thus, selecting $a$ involves balancing these two performance metrics: a small $a$ minimizes the total queue length for all patients, while a large $a$ prioritizes reducing the queue length for severe patients.

We show an example of determining the parameter $a$ as follows. 
We evaluate 75 parameter sets, using combinations of $\alpha \in \{0.05, 0.1, \ldots, 0.25\}$, $\beta \in \{0.8, 0.9, 1\}$, and $\gamma \in \{0.1, 0.2, \ldots, 0.5\}$. For each set $(\alpha, \beta, \gamma)$, we select the better policy between the shortest-first and longest-first rules for various values of $a \in \{0, 0.1, \ldots, 1\}$. Figure \ref{fig:parameter_tune} shows the average queue length for types 4 and 5 patients versus the average queue length for all patients. This tradeoff curve can help hospital managers compare performance measures under different parameters and choose the most suitable one. 
}


\section{Proof} \label{app:proof}

\subsection{Proof of Lemma 1 and Lemma 2} 

\begin{proof} [Proof of Lemma 1]   
We fix a policy $\pi$. If the policy $\pi$ does not depend on $a$, the term ${X}^\pi_{r,i,ns}(t)$ remains unchanged when we vary $a$.
Therefore, we have 
\begin{equation*}
  \begin{aligned}
J_i^\pi &= E \bigg[ \sum_{t=1}^T\sum_{r=1}^R h(r,a)   { X}^\pi_{r,i,ns}  (t)  \bigg]  \\
&= \sum_{r=1}^R h(r,a)  E \bigg[ \sum_{t=1}^T   { X}^\pi_{r,i,ns}  (t) \bigg],   \\
  \end{aligned}
\end{equation*}
where $E \bigg[ \sum_{t=1}^T   { X}^\pi_{r,i,ns}  (t) \bigg]$ remains fixed under the given policy $\pi$ when we vary $a$. Since the function $h(r,a) = r^a$ is increasing convex in $a$ and ${ X}^\pi_{r,i,ns}  (t)\ge 0 $, we have that $J_i^\pi$ is increasing convex in $a$. 
The results follow.   
\end{proof}   
In the interest of space, the complete proof for Lemma 2 is left to the online companion.

\subsection{Proof of Theorem 3} 
\begin{proof}  [Proof of Theorem 3]  
When $j=i$, it is easy to see that $s_1$ and $s_2$ are equivalent, and thus the result follows immediately. 
Now, we consider the case when $j>i$. 
 
Let $W_{m, 1,  ns, l }$ be the waiting time (excluding the service time in the queue) of patient $1$ when they enter the nurse queue for the $l$th time in system $s_m$; similar definition of $W_{m, 2,  ns, l }$ for patient 2. It is easy to see that 
 \begin{equation}\label{eq:w12ns_l}
 \begin{aligned}
   W_{1, 1,  ns, l}   &=  W_{2, 2,  ns, l} , \quad  l \in[i] ,   \\ 
    W_{1, 2,  ns, l}  &=    W_{2, 1,  ns, l}, \quad  l \in[i] , \\  
\end{aligned} 
\end{equation}  
using the \emph{swapping argument} and the equal service time used in the proof of Lemma 2. 
Further, we can show that when patient 2 in $s_1$ enters the system for the $k$th time ($i+1 \le k\le j$), patient 1 in $s_1$ has left the system already. Thus, $W_{1, 2,  ns, k  } = 0,  \quad i+1\le k \le j$ and 
  \begin{equation} \label{eq:w12ns_l1j}
 \begin{aligned}
{\sum_{l=1  }^{j  }   W_{1, 2,  ns, l }    } &=  {\sum_{l=1  }^{i  }   W_{1, 2,  ns, l }    },  \quad i+1\le k \le j .  \\ 
\end{aligned} 
 \end{equation}

Based on the specification of the holding costs incurred during the service times and those incurred during the waiting time at the nurse queue, as well as using the equal service time assumption, we can get 
\begin{equation}
 \begin{aligned}   
 c_{2}(a) -  c_{1}(a)  &=    { {\sum_{l=1  }^{i  } (i-l+1 )^{a}   W_{2, 1,  ns, l }    }     }  - {\sum_{l=1  }^{i  } (i-l+1 )^{a}  W_{1, 1,  ns, l }  
 + \sum_{l=1  }^{j } (j-l+1 )^{a}  W_{2, 2,  ns, l }  - {\sum_{l=1  }^{j  } (j-l+1 )^{a}   W_{1, 2,  ns, l }    }     } .
  \end{aligned}   
\end{equation}   

When $a=0$, we have 
\begin{equation} 
 \begin{aligned}   
 c_{2}(0) -  c_{1}(0)  &=    { {\sum_{l=1  }^{i  }    W_{2, 1,  ns, l }    }     }  - {\sum_{l=1  }^{i  }   W_{1, 1,  ns, l } + \sum_{l=1  }^{j }   W_{2, 2,  ns, l }   - {\sum_{l=1  }^{j  }   W_{1, 2,  ns, l }    }     }   \\ 
  &=    { {\sum_{l=1  }^{i  }    W_{2, 1,  ns, l }    }     }  - {\sum_{l=1  }^{i  }   W_{1, 1,  ns, l } + \sum_{l=1  }^{j }   W_{2, 2,  ns, l }   - {\sum_{l=1  }^{i  }   W_{1, 2,  ns, l }    }     }   \\ 
 &= {\sum_{l=1  }^{j  }  W_{2, 2,  ns, l }    }       {  - {\sum_{l=1  }^{i  }  W_{1, 1,  ns, l }    }     } \\ 
  &= {\sum_{l=1  }^{i }  W_{1, 1,  ns, l }    } + {\sum_{l=i+1  }^{j  }  W_{2, 2,  ns, l }    }       {  - {\sum_{l=1  }^{i  }  W_{1, 1,  ns, l }    }     } \\ 
  &={\sum_{k=i+1  }^{j  }  W_{2, 2,  ns, l }    } \\   
& \ge 0 ,
  \end{aligned}   
\end{equation} 
where the second equality is based on Equation \ref{eq:w12ns_l1j}, 
the third equality is based on Equation \ref{eq:w12ns_l}. 
It implies that the cost under $s_1$ is smaller than that under $s_2$. 
Hence, it is better to choose the shortest-first rule ($s_1$)  when $a=0$. 

When $a=1$, we have 
\begingroup 
\allowdisplaybreaks
 \begin{align*}   
 c_{2}(1) -  c_{1}(1)  &=  {  {\sum_{l=1  }^{i  } (i-l+1 )   W_{2, 1,  ns, l }    }     }  - {\sum_{l=1  }^{i  } (i-l+1 )  W_{1, 1,  ns, l }
 + \sum_{l=1  }^{j } (j-l+1 )  W_{2, 2,  ns, l } - {\sum_{l=1  }^{j  } (j-l+1 )   W_{1, 2,  ns, l }    }     } \\
 &=  {{\sum_{l=1  }^{i  } (i-l+1 )   W_{1, 2,  ns, l }    }     }  - {\sum_{l=1  }^{i  } (i-l+1 )  W_{2, 2,  ns, l } + \sum_{l=1  }^{j } (j-l+1 )  W_{2, 2,  ns, l }  - {\sum_{l=1  }^{i  } (j-l+1 )   W_{1, 2,  ns, l }    }     } \\
  &= - (j-i) \sum_{l=1  }^{i }  W_{1, 2,  ns, l }  - {\sum_{l=1  }^{i  } (i-l+1 )  W_{2, 2,  ns, l } +  \sum_{l=1  }^{i } (j-l+1 )  W_{2, 2,  ns, l } + \sum_{l=i+1  }^{j } (j-l+1 )  W_{2, 2,  ns, l }   } \\
    &= - (j-i) \sum_{l=1  }^{i } W_{1, 2,  ns, l }  +{ (j-i) \sum_{l=1  }^{i  }  W_{2, 2,  ns, l } + \sum_{l=i+1  }^{j } (j-l+1 )  W_{2, 2,  ns, l }   } \\
        &= - (j-i) \sum_{l=1  }^{j  } W_{1, 2,  ns, l }  +{ (j-i) \sum_{l=1  }^{i  }  W_{2, 2,  ns, l } + \sum_{l=i+1  }^{j } (j-l+1 )  W_{2, 2,  ns, l }   } \\ 
        &\le  - (j-i) \sum_{l=1  }^{j  } W_{1, 2,  ns, l }  +{ (j-i) \sum_{l=1  }^{i  }  W_{2, 2,  ns, l } + (j-i ) \sum_{l=i+1  }^{j }  W_{2, 2,  ns, l }   } \\ 
        &=  - (j-i) \sum_{l=1  }^{j  } W_{1, 2,  ns, l }  +{ (j-i) \sum_{l=1  }^{ j }  W_{2, 2,  ns, l }  } \\   
    &=   (j-i) \bigg(\sum_{l=1  }^{j }  W_{2, 2,  ns, l } - \sum_{l=1  }^{j }  W_{1, 2,  ns, l } \bigg)   ,
  \end{align*}   
   \endgroup  
  
where the second equality is based on Equation \ref{eq:w12ns_l}, the third equality results from combining $\sum_{l=1  }^{i  } (i-l+1 )   W_{1, 2,  ns, l } $ with $-\sum_{l=1  }^{i  } (j-l+1 )   W_{1, 2,  ns, l }$ into $- (j-i) \sum_{l=1  }^{i }  W_{1, 2,  ns, l }$, and splitting $\sum_{l=1  }^{j } (j-l+1 )  W_{2, 2,  ns, l }$ into $\sum_{l=1  }^{i } (j-l+1 )  W_{2, 2,  ns, l }$ and $\sum_{l=i+1  }^{j } (j-l+1 )  W_{2, 2,  ns, l }$, the fourth equality arises by combining $ - \sum_{l=1  }^{i  } (i-l+1 )  W_{2, 2,  ns, l } $ with $\sum_{l=1  }^{i } (j-l+1 )  W_{2, 2,  ns, l }$ into $(j-i) \sum_{l=1  }^{i  }  W_{2, 2,  ns, l } $, the fifth equality is based on Equation \ref{eq:w12ns_l1j}, the first inequality is based on $j-l+1\le j-i$ for $i+1\le l \le j$, and the sixth equality follows by combining  $(j-i) \sum_{l=1  }^{i  }  W_{2, 2,  ns, l }$ with $(j-i ) \sum_{l=i+1  }^{j }  W_{2, 2,  ns, l } $  into  $(j-i) \sum_{l=1  }^{j  }  W_{2, 2,  ns, l }$.


We know that patient 2 has the same time in service in $s_1$ and $s_2$.  
In Lemma \ref{lemma:discharge_sequence}, we have shown that 
\begin{equation}
\begin{aligned}
\sum_{l=1  }^{j }   W_{2, 2,  ns, l } +    {\sum_{l=1  }^{j  }  S_{2, 2,  ns, l}  }   = D_{2,  2,  ns, j   }    \le D_{1,  2,  ns, j   } = \sum_{l=1  }^{j }  W_{1, 2,  ns, l }   +  {\sum_{l=1  }^{j  }  S_{1, 2,  ns, l}  }    . 
\end{aligned} 
\end{equation} 
Thus, we have 
\begin{equation}
\begin{aligned}
\sum_{l=1  }^{j }  W_{2, 2,  ns, l } \le \sum_{l=1  }^{j }  W_{1, 2,  ns, l }    
\end{aligned} 
\end{equation} 
since $ {\sum_{l=1  }^{j  }  S_{2, 2,  ns, l}  } = {\sum_{l=1  }^{j  }  S_{1, 2,  ns, l}  }  $ based on the equal service time assumption. 

So, 
\begin{equation}
 \begin{aligned}   
 c_{2}(a) -  c_{1}(a) \le 0. 
  \end{aligned}   
\end{equation}   
It implies that the cost under $s_2$ is smaller than that under $s_1$. 
Hence, it is better to choose the long-first rule ($s_2$) when $a=1$. 

Furthermore, similar to the argument in Lemma \ref{lemma:convexity}, it is easy to show both $c_{2}(a) $ and $c_{1}(a)$ are increasing convex in $a$. Then there exists a unique value of $a$ at which $c_{2}(a) = c_{1}(a)$.  
This establishes the existence of a threshold value of $a=\hat{a}\in[0,1]$, such that 
the total cost incurred in $s_1$ is smaller than that under $s_2$ when $a\le \hat a$, and the total cost incurred in $s_2$ is smaller than that under $s_1$ when $a> \hat a$. 
 Thus, the result follows.   
 
\end{proof} 


\newpage 
\section*{Online companion}  

\noindent\textbf{Proof sketch.}  
The proof uses an induction argument. 
First, we use induction to show that $D_{1, 1, ns, \ell} \leq D_{1, 2, ns, \ell}$ for all $1\le \ell \le i$. This implies that patient 1 will be discharged earlier in system $s_1$ than patient 2 in system $s_1$ after finishing $i$ times of visit to the nurse queue. As an auxiliary result, we also show that the departure times from the content state have the same structure, i.e., $D_{1, 1, cs, \ell} \leq D_{1, 2, cs, \ell}$ for all $1\le \ell \le i$.   

{Second, leveraging this proved departure time comparison, we compare the discharge sequences of customers in systems $s_1$ and $s_2$. Here, we use the important assumption for the service time, i.e., the two patients have the same service times when they enter the needy (content) state for the first $l$th time ($1\le l \le i$) in systems $s_1$ and $s_2$, respectively. This allows us to \textit{swap} the patients in the two systems to compare their service times. In other words, during the first $i$ times of needy service and content service, patient 1 in system $s_1$ (or swapped to system $s_2$) shares an identical (cumulative) service time with patient 2 in system $s_2$ (or swapped to system $s_1$). Subsequently, for the subsequent $(j - i)$ times of needy service and content service, patient 2 in system $s_1$ and patient 2 in system $s_2$ share the same  (cumulative) service time.  As a result, the main difference in the total length of stay in the two systems is caused by delay at the nurse queue under different prioritization. In the proof we split the duration between the beginning of service and when both patients are discharged at both systems into two phases. That is, the interval starting from the initiation of service until the discharge of patient 1 in system s1 is designated as ``phase~1,''  Subsequent to phase 1, the interval until the discharge of both patients in both systems is termed ``phase~2.'' We then focus on the discharge time comparison in phase 2 to complete the proof.}

\noindent\textbf{Assumption.}  Under the deterministic service time setup, we assume that the service times for the two patients within systems $s_1$ and $s_2$, respectively, have the same service times when they enter the needy (content) state for the first $l$th time ($1\le l \le i$).   
Then we have that during the initial $i$ times of needy service and content service, patient 1 in system $s_1$ ($s_2$) shares an identical service time with patient 2 in system $s_2$ ($s_1$). Subsequently, for the subsequent $(j - i)$ times of needy service and content service, patient 2 in system $s_1$ and patient 2 in system $s_2$ share the same service time.   
Then we have  
\begin{equation}\label{eqn:assumption}  
\begin{aligned}
 S_{1, 1,  ns, l}   &=    S_{2, 2,  ns, l} =     S_{1, 2,  ns, l}  =    S_{2, 1,  ns, l}  , \quad  l \in[i] ,   \\ 
  S_{1, 1,  cs, l}   &=    S_{2, 2,  cs, l}   = S_{1, 2,  cs, l}   =    S_{2, 1,  cs, l}, \quad  l \in[i] ,\\   
  S_{1, 2,  ns, l + i }   &=    S_{2, 2,  ns, l +i } ,  \quad l \in[ j-i ] ,   \\ 
  S_{1, 2,  cs, l+i }   &=    S_{2, 2,  cs, l+i }  , \quad  l \in[j-i ]  .
 \end{aligned} 
 \end{equation}


\begin{proof}    [Proof of Lemma 2]  
Recall that we use $D_{m, 1,  ns, l }$ ($D_{m, 1,  cs, l }$) to denote the departure time when patient $1$ departs the needy (content) state for the $l$th time in system $s_m$; similarly for patient $2$. Patient 1 requires $i$ stage of nursing care, and patient 2 requires $j$ stage of nursing care, with $i\leq j$. We first use the induction method to show for system $s_1$ that 
\begin{equation} 
\begin{aligned}  
D_{1,  1,  ns, l }  & \le D_{1,  2,  ns, l }  ,\\ 
D_{1,  1,  cs, l }  & \le D_{1,  2,  cs, l }   . 
\end{aligned}
\end{equation} 
 for $l \in [i]$.    


\noindent\textbf{Base case.} 

For the base case of $l=1$, that is, when the patients 1 and 2 enter needy and content states for the first time, we have 
\begin{align}  
D_{1, 1,  ns, 1 } &= S_{1, 1,  ns, 1 }  ,  \\ 
D_{1, 2,  ns, 1 }    &= S_{1, 1,  ns, 1 }   + S_{1, 2,  ns, 1 }  , \\ 
D_{1,  1,  cs, 1 }  &= D_{1,  1,  ns, 1 }  +S_{1, 1,  cs, 1}   = S_{1,  1,  ns, 1 }  +S_{1, 1,  cs, 1} , \\
D_{1,  2,  cs, 1 }  &= D_{1,  2,  ns, 1 }  +S_{1, 2,  cs, 1}   =  S_{1, 1,  ns, 1 }   + S_{1, 2,  ns, 1 } +S_{1, 2,  cs, 1}    . 
\end{align}

It is easy to see that 
\begin{equation}
\begin{aligned} 
D_{1,  1,  ns, 1 }  &= S_{1,  1,  ns, 1 }  \\
     &\le  S_{1,  1,  ns, 1 }  +S_{1, 2,  ns, 1}  =D_{1,  2,  ns, 1 }   ,
\end{aligned}
\end{equation} 
and 
\begin{equation}
\begin{aligned} 
D_{1,  1,  cs, 1 }  &= S_{1,  1,  ns, 1 }  +S_{1, 1,  cs, 1}  \\
     &= S_{1,  1,  ns, 1 }  +S_{1, 2,  cs, 1}  \\
   &\le S_{1, 1,  ns, 1 }   + S_{1, 2,  ns, 1 } +S_{1, 2,  cs, 1}   =D_{1,  2,  cs, 1 }  .
\end{aligned}
\end{equation} 
Here, we used the assumption that the service time for patient 1 and patient 2 are the same when they enter the content state the first time.

Under the given priority, patient 1 finishes the service in the needy state first, and depart from content state (i.e., returns to the needy state) first as well in system $s_1$. Next, patient 1 enters the nurse queue for the second time before being discharged. 
So when $l=2$, we have  
\begin{align} 
D_{1, 1,  ns, 2 } &=  \max\{ D_{1,  1,  cs, 1 } , D_{1, 2,  ns, 1 }   \}    + S_{1, 1,  ns, 2 }  \label{eq:l2ns_1}   , 
\end{align} 
where the equality holds since patient 1 will be serviced by the nurse at time $\max\{D_{1, 1, cs, 1}, D_{1, 2, ns, 1}\}$. This timing accounts for two scenarios: (i) the needy queue is idle, (ii) patient 2 is still in the nurse queue when patient 1 transitions into the needy state for the second time.

Similarly, we have 
\begin{equation}  \label{eq:l2ns_2} 
\begin{aligned} 
D_{1, 2,  ns, 2 }  &=
   \max\{ D_{1,  2,  cs, 1 } , D_{1, 1,  cs, l } + S_{1, 1,  ns, l+1 }  \}   + S_{1, 2,  ns, 2}  \\
   &=
   \max\{ D_{1,  2,  cs, 1 } , D_{1, 1,  ns, l+1 }  \}   + S_{1, 2,  ns, 2}  \\  
\end{aligned}
\end{equation} 
if $ D_{1, 1,  cs, l } < D_{1,  2,  cs, 1 } \le D_{1, 1,  cs, l+1 } $ ($1\le l\le i-1$). 
Here,  
$ D_{1, 1,  cs, l } < D_{1,  2,  cs, 1 } \le D_{1, 1,  cs, l+1 } $ corresponds to the scenario where patient 2 leaves the content state for the first time after patient 1 has left the content state for the $l$-th time but before patient 1 leaves the content state for the $(l+1)$-th time ($l\le i-1$). 
In this scenario, patient 2 will be serviced by the nurse at time $\max\{D_{1, 2, cs, 1}, D_{1, 1, ns, l+1}\}$. This timing accounts for two scenarios: (i) the needy queue is idle, (ii) patient 1 is still in the nurse queue when patient 2 transitions into the needy state for the second time.



Then we have 
\begin{align} 
D_{1,  1,  cs, 2}  &= D_{1,  1,  ns, 2 }  +S_{1, 1,  cs, 2}  , \label{eq:l2ns_3} \\
D_{1,  2,  cs, 2 }  &= D_{1,  2,  ns, 2 }  +S_{1, 2,  cs, 2}   .  \label{eq:l2ns_4}   
\end{align} 
Then 
\begin{equation}\label{eq:l2ns_2_Dns2}
\begin{aligned} 
D_{1,  1,  ns, 2 }  &=   \max\{ D_{1,  1,  cs, 1 } , D_{1, 2,  ns, 1 } \}    + S_{1, 1,  ns, 2 }     \\  
&=   \max\{ D_{1,  1,  cs, 1 } , D_{1, 2,  ns, 1 }\}    + S_{1, 2,  ns, 2 }     \\  
&\le  D_{1,  2,  ns, 2 } , 
\end{aligned}
\end{equation}  
where the first equality is based on Equation~\ref{eq:l2ns_1}, the second equality is based on the assumption that $S_{1, 1,  ns, 2 }   = S_{1, 2,  ns, 2 } $, the first inequality is based on the induction conclusion $D_{1,  1,  cs, 1 } \le D_{1,  2,  cs, 1 } $ and $D_{1, 2,  ns, 1 }  \leq D_{1, 1,  ns, 2 } \leq D_{1, 1,  ns, 3} \le \cdots \le  D_{1, 1,  ns, i}  $.   


Similarly, we have 
\begin{eqnarray} 
\begin{aligned} 
D_{1,  1,  cs, 2 }  &= D_{1,  1,  ns, 2 }  + S_{1, 1,  cs, 2}  \\
     &= D_{1,  1,  ns, 2 }   + S_{1, 2,  cs, 2}  \\
   &\le   D_{1, 2,  ns, 2}   +  S_{1, 2,  cs, 2}       \\ 
   &=D_{1,  2,  cs, 1 }  , 
\end{aligned}
\end{eqnarray} 
where the first equality is based on Equation~\ref{eq:l2ns_3}, the second equality is based on the the assumption that $S_{1, 1,  cs, 2}  =S_{1, 2,  cs, 2} $,  the first inequality is based on Equation \ref{eq:l2ns_2_Dns2}, and the last equality is based on Equation~\ref{eq:l2ns_4}.

\noindent\textbf{Inductive step. } 

We assume it holds for $l=l'$, then we have 
\begin{equation} 
\begin{aligned}  
D_{1,  1,  ns, l' }  & \le D_{1,  2,  ns, l' }  ,\\ 
D_{1,  1,  cs, l' }  & \le D_{1,  2,  cs, l' }   . 
\end{aligned}
\end{equation} 
We have   
\begin{equation}  \label{eq:llp1ns_1}
\begin{aligned} 
D_{1, 1,  ns, l'+1  }  &=    \max\{ D_{1,  1,  cs, l' } , D_{1, 2,  ns, l''+1 }  \}   + S_{1, 1,  ns, l'}   ,  \\ 
\end{aligned}
\end{equation} 
if $ D_{1, 2,  cs, l'' } < D_{1,  1,  cs, l' } \le D_{1, 2,  cs, l''+1 } $ ($0\le l'' \le l'-1$), where $ D_{1, 2,  cs, 0 }=0$.   
Here,  
$ D_{1, 2,  cs, l'' } < D_{1,  1,  cs, l' } \le D_{1, 2,  cs, l''+1 } $ corresponds to the scenario where patient 1 leaves the content state for the $l'$th time after patient 2 has left the content state for the $l''$th time but before patient 2 leaves the content state for the $(l''+1)$th time ($0\le l''\le l'-1$). 
In this scenario, patient 1 will be serviced by the nurse at time $\max\{D_{1, 1, cs, l'}, D_{1, 2, ns, l''+1}\}$. This timing accounts for two scenarios: (i) the needy queue is idle, (ii) patient 2 is still in the nurse queue when patient 1 transitions into the needy state for the $ (l'+1)$th time.



Similarly, we have 
\begin{equation}  \label{eq:llp1ns_2}
\begin{aligned} 
D_{1, 2,  ns, l'+1 }  &=    \max\{ D_{1,  2,  cs, l' } , D_{1, 1,  cs, l'' } + S_{1, 1,  ns, l''+1 }  \}   + S_{1, 2,  ns, l'} \\ 
&=    \max\{ D_{1,  2,  cs, l' } , D_{1, 1,  ns, l''+1 }  \}   + S_{1, 2,  ns, l'} \\ 
\end{aligned}
\end{equation} 
if $D_{1, 1,  cs, l'' } < D_{1,  2,  cs, l' } \le D_{1, 1,  cs, l''+1 }  $ ($l'\le l''\le i-1$). 
Here,  
$ D_{1, 1,  cs, l'' } < D_{1,  2,  cs, l' } \le D_{1, 1,  cs, l''+1 } $ corresponds to the scenario where patient 2 leaves the content state for the $l'$th time after patient 1 has left the content state for the $l''$th time but before patient 1 leaves the content state for the $(l''+1)$th time ($l'\le l''\le i-1$). 
In this scenario, patient 2 will be serviced by the nurse at time $\max\{D_{1, 2, cs, l'}, D_{1, 1, ns, l''+1}\}$. This timing accounts for two scenarios: (i) the needy queue is idle, (ii) patient 1 is still in the nurse queue when patient 2 transitions into the needy state for the $ (l'+1)$th time. 

Then we have 
\begin{align} 
D_{1,  1,  cs, l'+1}  &= D_{1,  1,  ns, l'+1 }  +S_{1, 1,  cs, l'+1}  , \label{eq:llp1ns_3} \\
D_{1,  2,  cs, l'+1}  &= D_{1,  2,  ns, l'+1 }  +S_{1, 2,  cs, l'+1}   .  \label{eq:llp1ns_4}   
\end{align} 
Then 
\begin{equation}\label{eq:llp1ns_2_Dns3}
\begin{aligned} 
D_{1,  1,  ns, l'+1 }  &= \max\{ D_{1,  1,  cs, l' } , D_{1, 2,  ns, l''+1 }  \}   + S_{1, 1,  ns, l'+1} &\mbox{if $ D_{1, 2,  cs, l'' } < D_{1,  1,  cs, l' } \le D_{1, 2,  cs, l''+1 } $, $0\le l'' \le l'-1$} \\ 
&= \max\{ D_{1,  1,  cs, l' } , D_{1, 2,  ns, l''+1 }  \}   + S_{1, 2,  ns, l'+1} &\mbox{if $ D_{1, 2,  cs, l'' } < D_{1,  1,  cs, l' } \le D_{1, 2,  cs, l''+1 } $, $0\le l'' \le l'-1$} \\ 
&\le  D_{1,  2,  ns,  l'+1 } .  
\end{aligned}
\end{equation}  

where the first equality is based on Equation~\ref{eq:llp1ns_1}, the second equality is based on the assumption that $S_{1, 1,  ns,  l'+1 }   = S_{1, 2,  ns,  l'+1 } $, the first inequality is based on the induction conclusion $D_{1,  1,  cs,  l'} \le D_{1,  2,  cs,  l' } $, and  $D_{1, 2,  ns, l'' }  \leq D_{1, 1,  ns, l' } \leq D_{1, 1,  ns, l'+1} \le \cdots \le  D_{1, 1,  ns, i}  $.

Similarly, we have 
\begin{eqnarray} 
\begin{aligned} 
D_{1,  1,  cs,  l'+1 }  &= D_{1,  1,  ns,  l'+1 }  + S_{1, 1,  cs,  l'+1}  \\
     &= D_{1,  1,  ns,  l'+1 }   + S_{1, 2,  cs,  l'+1}  \\ 
   &\le   D_{1, 2,  ns,  l'+1}   +  S_{1, 2,  cs,  l'+1}       \\ 
   &=D_{1,  2,  cs,  l'+1 }  ,  
\end{aligned}
\end{eqnarray} 
where the first equality is based on Equation~\ref{eq:llp1ns_3}, the second equality is based on the assumption that $S_{1, 1,  cs,  l'+1 }   = S_{1, 2,  cs,  l'+1 } $, the first inequality is based on Equation \ref{eq:llp1ns_2_Dns3}, and the last equality is based on Equation~\ref{eq:llp1ns_4}.

So, it holds for $l=l'+1$ as well. 
Thus, 
\begin{equation} 
\begin{aligned}  
D_{1,  1,  ns, l }  & \le D_{1,  2,  ns, l }  ,\\ 
D_{1,  1,  cs, l }  & \le D_{1,  2,  cs, l }   . 
\end{aligned}
\end{equation} 
 for $l \in [i]$.   

 Let $l=i$, we have 
$$D_{1,  1,  ns, i  }   \le D_{1,  2,  ns,  i }   .  $$

\noindent\textbf{Comparing two systems. } To complete the proof, it remains to compare the departure times in the two systems $s_1$ and $s_2$. 

 We know that the service times for the two patients within systems $s_1$ and $s_2$, respectively, have the same service times when they enter the needy (content) state for the first $l$th time ($1\le l \le i$).   
That is, during the initial $i$ times of needy service and content service, patient 1 in system $s_1$ ($s_2$) shares an identical service time with patient 2 in system $s_2$ ($s_1$).  
Then we have 
$$  D_{2,  2,  ns, i  } = D_{1,  1,  ns, i  }  \le  D_{1,  2,  ns,  i } = D_{2, 1,  ns,  i }   .  $$    

It is easy to see that 
$$D_{m,  n,  ns, i  }  \le   D_{m, n,  ns,  j  } , \quad m,n \in \{1,2\} $$    
since the patient needs to receive an additional $j-i$ rounds of needy and content service after they leave the needy state for the $i$th time. Then we have 
\begin{equation}\label{discharge:1}
\begin{aligned}
D_{1,  1,  ns, i  }           \le D_{1,  2,  ns, j   }  , \\  
D_{1,  1,  ns, i  }           \le D_{2,  2,  ns, j   }  , \\   
D_{1,  1,  ns, i  }           \le  D_{2,  1,  ns, i   }  . \\   
\end{aligned}
\end{equation}

We see that patient 1 in system $s_1$ is discharged first.   
We can also see that in system $s_1$, patient 1 always receives the service first, and then is followed by patient 2 until patient 1 is discharged. 
Similarly, in system  $s_2$, {since the service times in each stage are the same for the two patients, swapping the two patients, we know that} patient 2 always receives the service first in $s_2$, and then is followed by patient 1 until either of them is discharged.  
So, we can split the duration between the beginning of service and when both patients are discharged at both systems into two phases. 
Specifically, the interval starting from the initiation of service until the discharge of patient 1 in system $s_1$ is designated as ``phase 1.'' Subsequent to phase 1, the interval until the discharge of both patients in both systems is termed ``phase 2.'' At the end of phase 1, the remaining times of needy service for patient 2 in system $s_1$ is $j - k$ ($1\le k\le i-1$), and the remaining times of needy service for patient 1 in system $s_2$ is $i-k$, either in the content or needy state. 
Here, $k$ is the number of times the patient 2 in $s_1$ goes through the needy  queue in phase 1. 
Additionally, we know that patient 1 in $s_1$ and patient 2 in $s_2$ have the same departure times from the content and needy states every time they exit these states in phase 1. As a result, once patient 1 in $s_1$ exits the system, patient 2 in $s_2$ simultaneously completes the $i$th needy service. Therefore, at the end of phase 1, the remaining number of needy service times for patient 2 in system $s_2$ is $j - i$,  who just enters the content state.

Figure \ref{example_lemma} provides a simple example to explain phases 1 and 2. For illustration, we assume that the service time in both the needy and content states for each patient is one period. Patient 1 has $r=2$ (with two stages), and patient 2 has $r=3$  (with three stages). In $s_1$, we serve patient 1 first, and in $s_2$, we prioritize patient 2. The figure show that in $s_1$, patient 1 initially receives the needy service in the first period, transitions to the content state with one stage remaining in the second period, rejoins the needy queue with one stage remaining in the third period, and ultimately departs from the system (getting discharged) by the end of the third period. On the other hand, patient 2 waits in the needy queue during the first period while patient 1 is being serviced. Patient 2 then receives service from the nurse starting from the second period (with three stages remaining), moves to the content state (with two stages remaining) in the third period, returns to the needy queue in the fourth period, enters the content state again (with one stage remaining) in the fifth period, rejoins the needy queue in the sixth period, and finally departs from the system (getting discharged) by the end of the sixth period. The transitions for patients 1 and 2 in $s_2$ are similar to that in system $s_1$. Comparing across both systems, we see that patient 1 in $s_1$ has the earliest departure time while patient 2 in $s_1$ has the latest departure time; the departure times of the two patients in $s_2$ are in between. The duration from the initiation of service (the beginning of the first period) to the discharge of patient 1 in $s_1$ (the end of the third period) is termed phase 1. The subsequent interval, from the end of phase 1 until the discharge of Patient 2 in $s_1$, is called phase 2.

\begin{figure}
         \centering
         \includegraphics[width=0.7\textwidth]{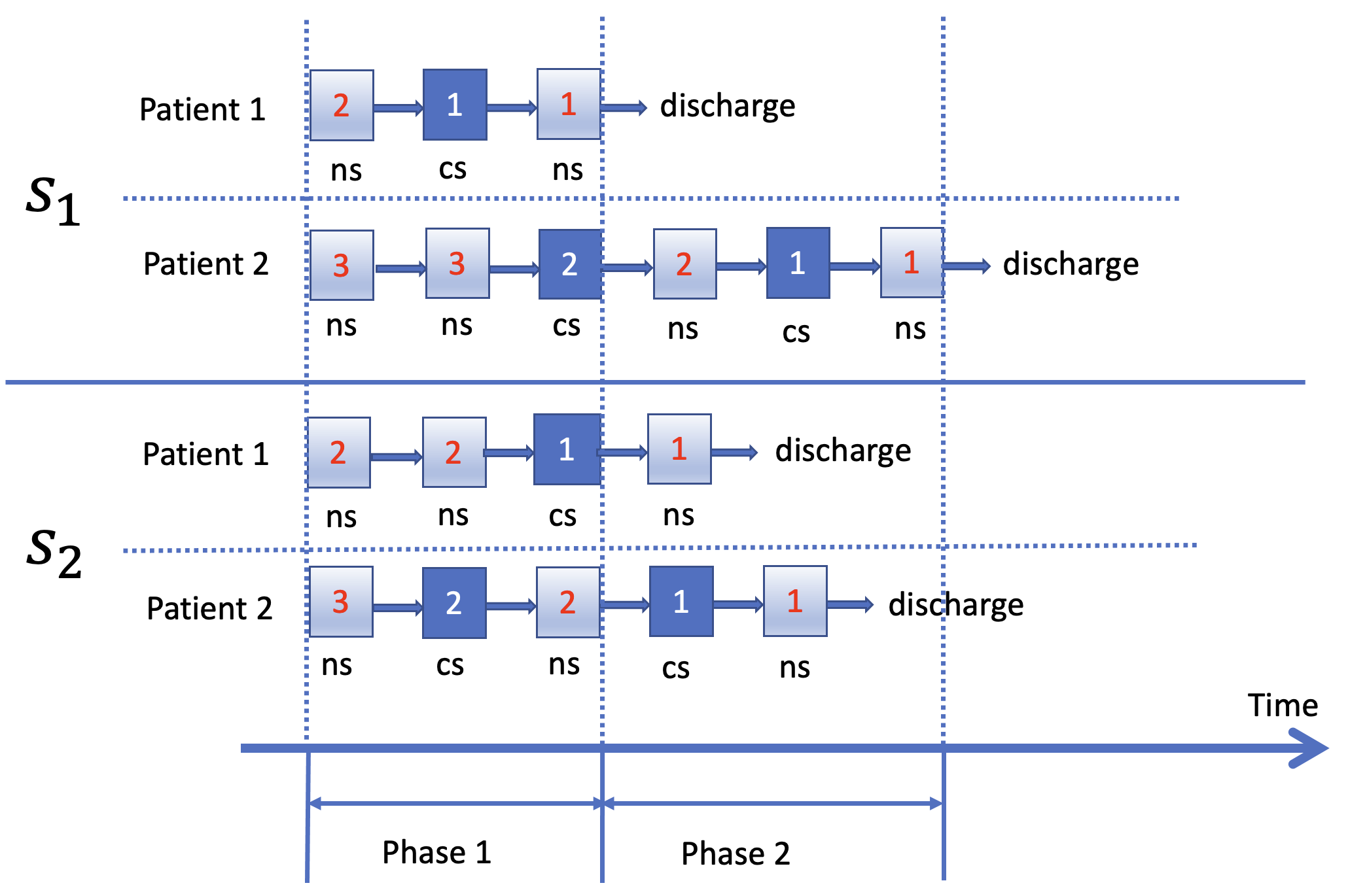}
        \caption{Two patient queueing system example}
        \label{example_lemma} 
\end{figure}



In phase 2, only patient 2 with type $j-k$ remains in the content or needy state in system $s_1$. 
Let $T_1$ represent the duration of phase 1. When patient 2 in $s_1$ is in the content state,  we have 
 \begin{equation}\label{eqn:discharge12nsj}   
\begin{aligned}
D_{1,  2,  ns, j }   =     T_1 +  \sum_{l=k+1  }^{ j}  S_{1,  2,  ns,  l }   +  \sum_{l=k+1   }^{j-1  }  S_{1,  2,  cs, l  }   +  S'_{1,  2,  cs, k  }  , 
\end{aligned}
\end{equation} 
where $S'_{1,  2,  cs, k   }  $ is the remaining service time in content state for patient 2 after phase 1.

Turning to system $s_2$,  patient 1 in $s_2$ is also in the content state when patient 2 in $s_1$ is in the content state. 
We let $S'_{2,  1,  cs, k  }  $ be the remaining service time in content state for patient 1 when patient 2 enters the content state for the $k$th time. 
Then we know $S'_{2,  1,  cs, k  } = S'_{1,  2,  cs, k   } $. 

 Then we have
 \begin{equation}\label{discharge:2}   
\begin{aligned}
D_{2,  1,  ns, i  } &  \le  T_1 +  \sum_{l=k+1 }^{ i}  S_{2,  1,  ns,  l }   +  \sum_{l=k+1 }^{ i-1}  S_{2,  1,  cs,  l }   +     S'_{2,  1,  cs,  k }  +  \sum_{l=i+1 }^{ j }    S_{2,  2,  ns,  l }     \\    
&  \le  T_1 +  \sum_{l=k+1 }^{ i}  S_{2,  1,  ns,  l }   +  \sum_{l=k+1 }^{ i-1}  S_{2,  1,  cs,  l }   +     S'_{2,  1,  cs,  k }  +  \sum_{l=i+1 }^{ j }    S_{2,  2,  ns,  l }     +   \sum_{l=i }^{ j -1}      S_{2,  2,  cs,  l }   \\    
&  =  T_1 +  \sum_{l=k+1 }^{ i}  S_{1,  2,  ns,  l }   +  \sum_{l=k+1 }^{ i-1}  S_{1,  2,  cs,  l }   +     S'_{1,  2,  cs,  k }  +  \sum_{l=i+1 }^{ j }    S_{1,  2,  ns,  l }     +   \sum_{l=i }^{ j -1}      S_{1,  2,  cs,  l }   \\   
&= T_1 +  \sum_{l=k+1  }^{ j}  S_{1,  2,  ns,  l }   +  \sum_{l= k+1  }^{j -1}  S_{1,  2,  cs, l  } +  S'_{1,  2,  cs, k }   \\ 
&=  D_{1,  2,  ns, j }  , 
\end{aligned}
\end{equation} 
where the first inequality holds since when patient 1 enters the needy state after phase 1, they will wait at most $\sum_{l=i+1 }^{ j }    S_{2,  2,  ns,  l }$, the second inequality holds since we add a nonnegative value $\sum_{l=i }^{ j -1}      S_{2,  2,  cs,  l }$, the first equality holds based on the model assumption in Equation \eqref{eqn:assumption}, and the last equality holds due to Equation~\eqref{eqn:discharge12nsj}. 


Similarly, for patient 2, we have 
  \begin{equation}\label{discharge:3}  
\begin{aligned}
D_{2,  2,  ns, j  }    &  \le T_1  +   \sum_{l=i+1 }^{ j }    S_{2,  2,  ns,  l }    +   \sum_{l=i }^{ j -1 }      S_{2,  2,  cs,  l } +  \sum_{l=k+1 }^{ i}  S_{2,  1,  ns,  l }    \\    
  &  \le T_1  +   \sum_{l=i+1 }^{ j }    S_{2,  2,  ns,  l }    +   \sum_{l=i }^{ j -1 }      S_{2,  2,  cs,  l } +  \sum_{l=k+1 }^{ i}  S_{2,  1,  ns,  l }  +  \sum_{l=k+1 }^{ i}  S_{2,  1,  cs,  l } + S'_{2,  1,  cs,  k }       \\    
    &  = T_1  +   \sum_{l=i+1 }^{ j }    S_{1,  2,  ns,  l }    +   \sum_{l=i }^{ j -1 }      S_{1,  2,  cs,  l } +  \sum_{l=k+1 }^{ i}  S_{1,  2,  ns,  l }  +  \sum_{l=k+1 }^{ i}  S_{1,  2,  cs,  l } + S'_{1,  2,  cs,  k }       \\    
&= T_1 +  \sum_{l=k+1  }^{ j}  S_{1,  2,  ns,  l }   +  \sum_{l= k+1  }^{j -1}  S_{1,  2,  cs, l  } +  S'_{1,  2,  cs, k }   \\ 
&=  D_{1,  2,  ns, j }  , 
\end{aligned}
\end{equation}  
where the first inequality holds since when patient 2 enters the needy state after phase 1, they will wait at most $\sum_{l=k+1 }^{ i}  S_{2,  1,  ns,  l }$, the second inequality holds we add a nonnegative value $\sum_{l=k+1 }^{ i}  S_{2,  1,  cs,  l } + S'_{2,  1,  cs,  k } $, the first equality holds based on the model assumption in Equation \eqref{eqn:assumption}, and the last equality holds due to Equation \eqref{eqn:discharge12nsj}. 

Similarly, in cases where both patient 2 in $s_1$ and patient 1 in $s_2$ are in the needy state, we can similarly show 
 $D_{2,  1,  ns, i  }   \le   D_{1,  2,  ns, j }  $ and $D_{2,  2,  ns, j  } \le   D_{1,  2,  ns, j }  $.

In summary, we have 
\begin{equation}
\begin{aligned}
 D_{1,  1,  ns, i  }      &\le D_{2,  1,  ns, i   }      \le  D_{1,  2,  ns, j   }   , \\ 
 D_{1,  1,  ns, i  }    & \le D_{2,  2,  ns, j   }       \le D_{1,  2,  ns, j   }    . 
\end{aligned} 
\end{equation} 

The result follows. 
\end{proof}

\end{document}